\documentclass[12pt,reqno]{amsart}

\usepackage{AKstyle}

\usepackage[final,
color]{showkeys}
\definecolor{refkey}{gray}{.85}
\definecolor{labelkey}{gray}{.85}

\numberwithin{equation}{section}

\usepackage[pagebackref=true, colorlinks]{hyperref}

\hypersetup{pdffitwindow=true,linkcolor=blue,citecolor=blue,urlcolor=blue,menucolor=blue}

\usepackage{comment}

\begin{document}

\title[Symmetry Types of  $\GL(3)$ families]{
On the $\GL(3)$ Kuznetsov Formula
with applications to 
%The
 Symmetry Types of 
families of $L$-functions
% of 
% $\GL(3)$, $\GL(6)$, $\GL(8)$, and $\GL(9)$ %Cuspidal, Symmetric Square, Adjoint, and Rankin-Selberg 
}
%{\tiny\it\rm very preliminary version, please do not distribute}
\author{Dorian Goldfeld}
\address{Math Department \\ Columbia University\\ New York, NY 10027}
\email{goldfeld@math.columbia.edu}

\author{Alex Kontorovich}
\address{Math Department \\ Yale University\\ New Haven, CT 06511}
\email{alex.kontorovich@yale.edu}

\thanks{D.G. is partially supported by  NSF grant DMS-1001036.}
\thanks{A.K. is partially supported by NSF grants  DMS-1209373, DMS-1064214 and DMS-1001252.}

\begin{abstract}
We present an explicit approach to the $\GL(3)$ Kuznetsov formula. 
As an application, for a restricted class of test functions, we obtain the low-lying zero densities for the following three families:
cuspidal $\GL(3)$ Maass forms $\phi$, the  symmetric square family $\sym^{2}\phi$ on $\GL(6)$, and 
the adjoint family $\Ad \phi$ on $\GL(8)$. 
Hence we can identify their symmetry types; they
 are: unitary, unitary, and symplectic, respectively.
\end{abstract}
\date{\today}
\maketitle
\tableofcontents

\section{Introduction}

\subsection{Symmetry Types}
\
\vskip 5pt
% Motivated by function field analogues,  
 In \cite{KatzSarnak1999},
 Katz and Sarnak 
 introduced the notion of symmetry type for a family of $L$-functions.
 Since then there has been a slew of activity  
regarding the following problem: given a family of $L$-functions,   
  determine its symmetry type.
 %Let $\mathcal F$ be a family of automorphic forms. To each $f \in \mathcal F$ we associate an $L$-function $$L(s,f) = \sum_{n=1}^\infty A _f(n) n^{-s}, \qquad (A_f(n) \in \Bbb C),$$
% which converges absolutely for $\Re(s)$ sufficiently large. By abuse of notation, we shall also refer to $\mathcal F$ as a family of $L$-functions.
A common approach
to this determination,
 as outlined in \cite{Sarnak2008}, is to 
 analyze the density of low-lying zeros in  the specified family, for test functions whose Fourier transforms have restricted support. 
 Such an analysis has been carried out  in many places, including \cite{IwaniecLuoSarnak2000, Royer2001, Gulolu2005, Young2006, HughesMiller2007, MillerEtAl2011}
 for $\GL(2)$ and \cite{DuenezMiller2006} for some $\GL(4)$ and $\GL(6)$ families; see also \cite{DuenezMiller2009}. 
In all cases in the literature (going beyond $GL(1)$ or abelian methods), the analysis involves a version of the $\GL(2)$ Petersson/Kuznetsov formula.

 The purpose of this paper is to carry out a similar analysis, for the first time using the $\GL(3)$ Kuznetsov formula.
 Assuming the generalized Riemann hypothesis (to interpret the low-lying zero densities) and the generalized Ramanujan conjectures (for ease of exposition), we will determine symmetry types for the families of 
 \begin{enumerate}
 \item
  cuspidal $\GL(3,\Z)$ automorphic forms $\phi$,
  \item
   the $\GL(6)$ symmetric square family $\sym^{2}\phi$, and
   \item the $\GL(8)$ adjoint family $\Ad\phi$.
   %, where if $\pi$ is the corresponding automorphic representation, then we have the isobaric factorization $\pi\times\tilde\pi=\Ad\pi\boxplus\bo$. (Here $\bo$ is the trivial representation corresponding to the classical Riemann zeta function).
   %, and lastly
%   \item the $\GL(9)$  Rankin-Selberg family $\pi\times\pi$. Of course this family factors generically as $\pi\times\pi=\sym^{2}\pi\boxplus\wedge^{2}\pi$, so the low-lying zero density is just the sum of the factors. On $\GL(3)$, the exterior square is simply the contragredient representation, 
%   so 
%the Rankin-Selberg square symmetry type is obtained immediately from knowledge of $(1)$ and $(2)$.
   %we do not study it separately. 
\end{enumerate}
We will show that the symmetry types are:
 \begin{enumerate}
 \item
unitary,
  \item
unitary,
   and 
   \item 
   symplectic,
%   \item
%unitary,
\end{enumerate}
respectively. 

The methods presented here are capable of wide generalization, in particular, it should be possible to determine the symmetry types of families associated to  pairs of cuspidal automorphic representations on $GL(n)$ for any $n\ge 2.$ We hope to return to this topic in a future publication.
\vskip 5pt
To state our results more precisely, we need some background.
\\

\subsection{Hecke-Maass forms}\
\vskip 5pt
Let $G=\GL(3,\R)$ with maximal compact $K=O(3)$ and center $Z=\R^{\times}$, 
let $\fh^{3}=G/(K\cdot Z)$ be the generalized upper half plane,
and take the lattice $\G:=\GL(3,\Z)$ in $G$. %\footnote{% 
%Note that we do not work adelically, because the Poincar\'e series which we use
% %to prove 
% in the $\GL(3)$ analog of  Kuznetsov's formula do not exist adelically.
%That said, the formula is a consequence of the relative  formula, and hence can be computed adelically, see e.g. \cite{Ye2000}. We prefer the % explicit 
%treatment with Poincar\'e series, since our analysis requires explicit estimates.}

The algebra of $G$-invariant differential operators  
acts
on  $\cH:=L^{2}(\G\bk\fh^{3})$.
The Hecke-Maass forms $\phi_{j}$  constitute an orthogonal (Hecke normalized) basis for
$$
\cH_{0}:=\bigoplus_{j=1}^\infty\C\phi_{j}\subset \cH
,
$$ 
where $\cH_{0}$ is the cuspidal
subspace in the 
Langlands spectral decomposition
%of, given by
 \cite[Prop. 10.13.1]{Goldfeld2006}
$$
\cH%L^{2}(\G\bk \fh^{3})
=
\C\bo\oplus
\cH_{0}
\oplus
\cH_{min}
\oplus
\cH_{max}
\oplus
\cH_{res}
.
$$
%see \S\ref{subs:spec} below. 
Here % $\cH_{0}$ is the cuspidal spectrum, and 
$\cH_{min},$ $\cH_{max},$ and  $\cH_{res}$ are, respectively, the spans of integrals of the minimal and maximal parabolic Eisenstein series, and the residual spectrum.  

%By the Duality Theorem
%\cite{GelfandGraevPS1966},
%the $G$-span of a Maass form $\phi$ gives rise to an irreducible representation $(\pi,V)$
%in the decomposition of the right-regular representation of $G$ acting on $L^{2}(\G\bk G)$.
%a principal series representation with purely imaginary 

Let the Hecke-Maass form $\phi_{j}$ have
spectral parameters $\nu^{(j)}:=(\nu_{1}^{(j)},\nu_{2}^{(j)})$. 
When discussing a fixed form $\phi$, we drop the superscripts $(j)$.
Our normalization\footnote{Note that our normalization differs from that used in \cite{Goldfeld2006} by  $\nu_{j}\mapsto1/3+\nu_{j}$.} 
is such that 
% the Ramanujan conjectures% at infinity
 %, assumed henceforth (for ease of exposition), imply that
for a tempered form,
  $\nu_{1}$ and $\nu_{2}$ are purely imaginary.

It is convenient to also introduce the spectral parameters
$$
\nu_{3}:=\nu_{1}+\nu_{2},
$$
and
$$
\ga_{1}:=\nu_{1}+\nu_{3},\quad
\ga_{2}:=-\nu_{1}+\nu_{2},\quad
\ga_{3}:=-\nu_{2}-\nu_{3}.
$$
Writing $\gl_{\phi}$ for the Laplace eigenvalue of $\phi$, we have
$$
\gl_{\phi}
=
1-3(\nu_{1}^{2}+\nu_{2}^{2}+\nu_{3}^{2})
=
1-(\ga_{1}^{2}+\ga_{2}^{2}+\ga_{3}^{2})
%=
%1
%+3(t_{1}^{2}+t_{1}t_{2}+t_{2}^{2})
.
$$
%where we have set
%$$
%\nu_{1}=i t_{1},
%\qquad
%\nu_{2}=i t_{2}.
%$$

Weyl's Law in this setting \cite{Miller2001} states that 
$$
\#\{\phi:\gl_\phi<T^{2}\}\sim c T^{5}
,
$$
as $T\to\infty$, for some constant $c > 0.$
\\

\subsection{The $\GL(3)$ Kuznetsov Formula}\
\vskip 5pt
We
will state and use the $\GL(3)$ Kuznetsov formula with some naturally occurring weights, defined as follows.
For $j = 1,2,3,\ldots,$ let 
$$
\cL_{j}:=
\underset{s=1}{\Res}\
L\big(s, \, \phi_{j}\times\tilde\phi_{j}\big)
$$
be the residue at the edge of the critical strip of the $L$-function attached to $\phi_{j}\times\tilde\phi_{j}$;
generically this is the value at $s=1$ of $L(s, \Ad \phi_{j})$.

We introduce an absolute constant 
$$
R\ge 10,
$$ 
which is needed for certain technical reasons, see the estimates in \S\ref{subsec:pEst}.
For $T\gg 1$, we define 
\be\label{eq:hTdef}
h_{T,R}(\nu)
:=
e^{\left(\ga_{1}^{2}+\ga_{2}^{2}+\ga_{3}^{2}\right)/T^{2}}\;
{\left(
\prod\limits_{1\le j\le 3} 
\G\left({2+R+3\nu_j\over4}\right)
\G\left({2+R-3\nu_j\over4}\right)
\right)^{2}
\over
\prod\limits_{1\le j\le 3} 
\G\left({1+3\nu_j\over2}\right)
\G\left({1-3\nu_j\over2}\right)
}
.
\ee
In fact, $R$ is needed to enable us later to pull  contours in certain integrals with respect to the $\nu_{j}$'s without passing through poles of the numerator in \eqref{eq:hTdef}.
Note that $h_{T,R}(\nu)>0$, and is essentially supported on $\gl<T^{2},$ or $|\nu_{1}|,|\nu_{2}|,|\nu_{3}|\ll T$. 
In this range, 
one sees from Stirling's formula 
that
if $\nu$ is tempered, then
%$h_{T,R}$ is asymptotic to a constant times
%$
%C_{\nu}^{R},
%$
%where we have
%set
\be\label{eq:hAsymp}
%C_{\nu}:=
h_{T,R}(\nu)
\sim 
c_{R}
\big[
(1+|\nu_{1}|)(1+|\nu_{2}|)(1+|\nu_{3}|)
\big]^{R}
,
\ee
  for some $c_{R}>0$.
The non-tempered forms constitute a zero density set \cite{Miller2001}.
 
 Let $A_j
(n_{1},n_{2})$ denote  the coefficients of $\phi_{j}$  
%according to Piatetski-Shapiro's
in the Fourier-Whittaker
 expansion, see \S\ref{subs:cusp}. 

\begin{thm}\label{thm:Kuz}
With the above notation and assuming the Ramanujan conjecture at the infinite place%
, we have the ``Weyl Law'', that for some $c>0$,
\be\label{eq:hTmain}
\sum_{j}{h_{T,R}(\nu^{(j)})\over \cL_{j}} 
\sim c\ T^{5+3R}.
\ee
Moreover for fixed $\vep > 0$, $R\ge 10$, $n_{1},n_{2},m_1,m_2\in\Z_{\ge1}$, and $T\gg 1$, we have
\be\label{eq:KuzAns}
\sum_{j}A_{j}(m_{1},m_{2}) \overline{A_j(n_1,n_2)}\,
%\overline{A_{j}(m_{1},m_{2})}
{h_{T,R}(\nu^{(j)})\over \cL_{j}} 
=
%\bo_{
%n_{1}=n_{2}=1\atop m_1=m_2=1
%n_{1}=m_{1}\atop n_{2}=m_{2}
%}
\begin{cases} \sum\limits_j {h_{T,R}(\nu^{(j)})\over \cL_{j}}+
 \cO_{R,\vep}
\Bigg(
T^{3+3R+\vep} \,
|m_1m_2n_{1}n_{2}|^2
\Bigg), & \text{if } {m_1=n_1,
\atop 
m_2=n_2,
}\\
\cO_{R,\vep}
\Bigg(
T^{3+3R+\vep}\,
|m_1m_2n_{1}n_{2}|^2
\Bigg), &\hskip -20pt {otherwise.}\end{cases}
\ee
\end{thm}

\begin{rmk}
In light of the asymptotic formula \eqref{eq:hAsymp}, the analytic weight $h_{T,R}$ can %surely 
be removed with a modicum of effort; 
we have chosen to leave the weight
for ease of exposition. 
The same is done for $\GL(2)$ in \cite{MillerEtAl2011}.
\end{rmk}

\begin{rmk}
The weight $\cL_{j}$ is more subtle; 
it is shown in \cite[(1.4)]{Blomer2011} that
$$
C_{\nu^{(j)}}^{-1}\ll \cL_{j} \ll_{\vep} C_{\nu^{(j)}}^{\vep},
$$
where
$$
C_{\nu}=(1+|\nu_{1}|)(1+|\nu_{2}|)(1+|\nu_{3}|)
.
$$
Moreover if one assumes the functorial transfer predicting $\phi\times\tilde\phi$ is automorphic on $\GL(9)$, then using
the non-existence of Siegel zeros for the corresponding $L$-function
\cite{HoffsteinRamakrishnan1995}, 
one can improve 
the lower bound above to $C_{\nu}^{-\vep}$.
With this assumption, the weight can be removed completely, as in \cite{Luo2001},
giving rise to a clean cut-off.
\end{rmk}

\begin{rmk}\label{rmk:imp}
We have not made any attempt to obtain the best possible error terms in \eqref{eq:KuzAns}. In particular, we have made no use of stationary phase, nor have we even invoked Deligne's bounds for Kloosterman sums (see e.g. \cite[Larsen's appendix]{BumpFriedbergGoldfeld1988}). We tried to present as simple a method as we 
could, keeping in mind the eventual goal of generalizing these techniques to $\GL(n)$ with $n \ge 2.$
\begin{comment}
The error terms in \eqref{eq:KuzAns} can surely be improved with a modicum of effort (indeed, as remarked above, for less explicit functions $h$, Blomer's method gives better errors). In particular, we have made no use of stationary phase, nor have we even invoked Deligne's bounds for  Kloosterman sums (see e.g. \cite[Larsen's appendix]{BumpFriedbergGoldfeld1988}). We tried to present as simple a  method as we could, keeping in mind the eventual goal of generalizing these techniques to $\GL(n)$.
\end{comment}
\end{rmk}

\begin{rmk}
A similar % but less precise 
result is obtained in \cite{Blomer2011}. %. Blomer's analysis
%with some key distinctions.
Blomer first chooses a test function on the geometric side, and then executes a delicate analysis to obtain implications on the spectral side. In our approach, we choose the test function on the spectral side first, making the asymptotic formula  \eqref{eq:hAsymp} 
immediately visible. In a private communication, Blomer has informed us that 
%by extending 
from
%refining
the methods in \cite{Blomer2011}, he can also obtain \eqref{eq:hTmain}  and \eqref{eq:KuzAns}  (with a better error term) for a range of test functions.
\begin{comment}
Blomer's analysis proceeds by cleverly picking a test function on the arithmetic side, and
executing
 a delicate analysis of its implications on the spectral side. For various choices of his test functions,  the spectral side (that is, the function $h$ in \eqref{eq:KuzAns}) can be made to isolate 
Maass forms with spectral parameters near a specified location. 
In fact his method gives better error terms than what we obtain in \eqref{eq:KuzAns}.
But he does not give asymptotic formulae for $h$ as we have in \eqref{eq:hAsymp}, since they are not needed in his applications; he evaluates his $h$ up to constant. In our application to low-lying zeros, the exact asymptotics are crucial, since a constant fluctuation could {\it a priori} influence the symmetry type. Hence we proceed instead by picking explicit test functions on the spectral side, as in \eqref{eq:hTdef}, and using the $\GL(3)$ analog of the Kontorovich-Lebedev inversion (see \S\ref{subs:KLt}), %as analyzed in \cite{GoldfeldKontorovich2011}, 
to give bounds on the arithmetic side.
That said, Blomer's method can also be made to give an asymptotic formula for his $h$. % with more work.
\end{comment}
\end{rmk}

\begin{comment}
\begin{rmk}\label{rmk:RamUse}
The Ramanujan conjectures are not used in any essential way in the proof of Theorem \ref{thm:Kuz}; just to simplify the exposition. Since we will be assuming GRH to interpret the low-lying zero densities, we may as well allow ourselves (and the reader) the convenience.
\end{rmk}
\end{comment}

\subsection{Low-Lying Zeros}\
\vskip 5pt
For a Hecke-Maass form $\phi$ on $\GL(3)$, let $\rho(\phi)$ be one of
$$
\rho(\phi)
=
\threecase
{\phi}{}
{\sym^{2}\phi}{}
{\Ad\phi,}{}
%{\phi\times\phi,}{}
$$
and let $L(s,\rho(\phi))$ be the corresponding $L$-function. Let $\alpha_1, \alpha_2, \alpha_3$ be the spectral parameters 
associated to $\phi$. If the Laplace eigenvalue $\lambda_\phi = 1- (\alpha_1^2+\alpha_2^2+\alpha_3^2)$ is sufficiently large, then we define the analytic conductor $c_{\rho(\phi)}$ of $\rho(\phi)$  as follows.
\be\label{eq:AnalyCond}
c_{\rho(\phi)}
=
\begin{cases}   
%\log\left( 
\pi^{-3}\cdot \underset {|\alpha_k| \,\ge \,\frac12} {\prod\limits_{1\le k\le 3}}\frac{|\alpha_k|}{2}
%\right)
, &\text{ if } \rho(\phi) = \phi,\\
& \\
 %\log\left( 
 \pi^{-5}\cdot \underset {|\alpha_j+\alpha_k| \,\ge \,\frac12} {\prod\limits_{1\le j\le k\le 3}}\frac{|\alpha_j+\alpha_k|}{2}
 %\right)
 , &\text{ if } \rho(\phi) = \sym^{2}\phi,\\
 & \\  
  %\log\left( 
  \pi^{-9}\cdot \underset {|\alpha_j-\alpha_k| \,\ge \,\frac12}{\prod\limits_{j=1}^3\prod\limits_{k=1}^3} \frac{|\alpha_j-\alpha_k|}{2}
  %\right)
  , & \text{ if }\rho(\phi) = \phi\times \bar\phi. \end{cases}
\ee

\begin{rmk}
Note that this is off by a constant from the more standard Iwaniec-Sarnak definition of ``conductor,'' for which see e.g. \cite[p. 95]{IwaniecKowalski}. The constants are crucial in our applications, see specifically \eqref{eq:ArhoPhi}, so we make our definition as above.
\end{rmk}

 We are interested in the weighted average value, denoted $C_\rho$, of the conductor $c_{\rho(\phi)}$  with respect to the weighting function $h_{T,R}$ defined in \eqref{eq:hTdef}. Then $C_\rho$ is defined by
 \begin{equation} \label{Arho}
\sum_j 
%\frac{
\log c_{\rho(\phi_j)}
%}{
{h_{T,R}(\nu^{(j)})\over \cL_{j}} \; \sim \; 
 %\left(
\log C_\rho \sum_{j}{h_{T,R}(\nu^{(j)})\over \cL_{j}} 
%\right)^{-1} 
, \qquad(T\to \infty),
\end{equation}and satisfies
$$
C_\rho
\asymp
\threecase
{T^{3}}{if $\rho(\phi)=\phi,$}
{T^{6}}{if $\rho(\phi)=\sym^{2}\phi,$}
{T^{6}}{if $\rho(\phi)=\Ad\phi$.}
%{T^{9}}{if $\rho(\phi)=\phi\times\phi$.}
$$
The weighted average value of the conductor in a family is introduced to normalize the low-lying zeros for comparison between the different families and the different matrix ensembles.

Let $\psi$ be an even test function of Schwartz class on $\R$ and define the low-lying zeros sum
$$
D(\rho(\phi); \,\psi):= \sum_{\g}\psi\left(\g \,{\log C_{\rho}\over 2\pi} \right),
$$
where $\g$ runs over the ordinates of nontrivial zeros of $L(s,\rho(\phi))$, counted with multiplicity.
To interpret this as capturing the low-lying zeros, we must assume GRH for the corresponding $L$-functions.
As $\psi$ has rapid decay, this sum localizes to those $\g$ which are within $1/\log C_{\rho}$ of the origin (corresponding to the central point $s=1/2$ of the $L$-function). 
 
% To compute the distribution of low-lying zeros, we must determine the limiting density  $W$ for which  $D(f;\psi)$, averaged over $f\in\cF_{R S}(T)$ with some spectral weights $w_f$, has a limit as $T\to\infty$, that is, 

\begin{theorem}\label{thm:main}
%Assume that the family $\mathcal F_{R S}$ of Rankin-Selberg L-functions  for $SL(3, \mathcal Z)$ satisfies Assumptions 1,2,3. If  
Assume the Fourier transform $\hat\psi$ of $\psi$ has   support in $(-\delta, \delta)$,
where
$$
\delta
=
\threecase
{4/15,}{if $\rho(\phi)=\phi,$}
{2/27,}{if $\rho(\phi)=\sym^{2}\phi,$}
{2/27,}{if $\rho(\phi)=\Ad\phi$.}
%{1/90}{if $\rho(\phi)=\phi\times\phi$.}
$$
Assume the Ramanujan conjectures, and GRH for the corresponding $L$-functions.
Then we have  the asymptotic formula
% \eqref{spw} holds. 
\be\label{spw}
{
1
\over
\sum_{j} {h_{T,R}(\nu^{(j)})\over \cL_{j}}
} \cdot
\sum_{j}D(\rho(\phi_{j}); \,\psi)\, {h_{T,R}(\nu^{(j)})\over \cL_{j}}
 \ =\int_\R \psi(x) W_{\rho(\phi)}(x)dx  + \mathcal O\left(\frac{\log\log T}{\log T}  \right) 
,
\ee
as $T\to\infty$, with
the limiting density function $W$ above 
%is 
given by
\be\label{Wis}
W_{\rho(\phi)}(x)
=
\threecase
{1,}{if $\rho(\phi)=\phi,$}
{1,} {if $\rho(\phi)=\sym^{2}\phi,$}
{1-{\sin(2\pi x)\over 2\pi x},}{if $\rho(\phi)=\Ad\phi$.}
%{1-{\sin(2\pi x)\over 2\pi x}}{if $\rho(\phi)=\phi\times\phi$.}
\ee
That is, the family $\rho(\phi)$  has symmetry type: unitary, unitary, %orthogonal, 
and symplectic, respectively.
\end{theorem}

\begin{remark}
The exterior square L-function on $GL(3)$ is the same as the contragredient L-function (see  \cite{BumpFriedberg1990, JacquetShalika1990, Kontorovich2010}). So the symmetry type for the exterior square family is unitary.
\end{remark}

\begin{rmk}
Note that \eqref{Wis} is consistent with  a recent conjecture by Shin and Templier \cite{ShinTemplier2012}.
\end{rmk}

\begin{rmk}
As in Remark \ref{rmk:imp}, the range of $\gd$ above can also be improved, and is intimately tied to the error terms in \eqref{eq:KuzAns}. 
%In fact, using the approach in \cite{Blomer2011}, Blomer has pointed out to us that the 
\end{rmk}

\begin{rmk}
The Ramanujan conjectures % at the infinite and finite primes 
are assumed to make the exposition of Theorem \ref{thm:main} as simple as possible. They can easily be removed by decreasing the size of $\delta$ in Theorem \ref{thm:main}.
\begin{comment}
As in Remark \ref{rmk:RamUse}, the Ramanujan conjectures are used here to simplify the argument. The most crucial place where this assumption is used is in dropping the terms with $\ell\ge3$ from \eqref{Dphipsi}; see just before \eqref{DphipsiRevised}. The Jacquet-Shalika bound \cite{JacquetShalika1981} or even the Luo-Rudnick-Sarnak bound \cite{LuoRudnickSarnak1999} is insufficient to avoid analyzing contributions from higher prime powers, which leads in the end to worse restrictions on the support $\gd$.
\end{comment}
\end{rmk}

\begin{comment}
By Parseval's Theorem,
$$
\int_\R\psi(x)W_{\rho(\phi)}(x)dx = \int_\R \hat\psi(y)\hat W_{\rho(\phi)}(y)dy,
$$
where $\hat\psi$ is the Fourier transform of $\psi$, and $\hat W_{\rho(\phi)}$ can be explicitly computed from \eqref{Wis} as a distribution:
$$
\hat W_{\rho(\phi)}(y)=\begin{cases} \delta_0(y), & {\rm if } \;\rho(\phi)=\phi,\\
& \\
\delta_0(y), & {\rm if } \;\rho(\phi)=\sym^{2}\phi,\\
&
\\

  \gd_0(y)-\scriptstyle{\threecase{\scriptstyle1/2,}{$\scriptstyle |y|<1,$}{\scriptstyle 1/4,}{$\scriptstyle |y|=\pm1,$}{\scriptstyle 0,}{ $\scriptstyle |y|>1,$.}} & {\rm if} \; \rho(\phi) =\Ad\phi.\end{cases}
$$
As the support of $\hat\psi$ will be restricted well inside $(-1,1)$, it follows that
\begin{equation}\int_\R \hat\psi(y)\hat W_{\rho(\phi)}(y)dy = \begin{cases} \int_\R \hat\psi(y)\gd_0(y)dy = \widehat\psi(0), & {\rm if} \; \rho(\phi) = \phi,\\
\int_\R \hat\psi(y)\gd_0(y)dy = \widehat\psi(0), & {\rm if} \; \rho(\phi) = \sym^2\phi\\
\int_\R \hat\psi(y)\gd_0(y)dy - \int_{supp(\hat\psi)}\hat\psi(y) \foh dy,  \\
\hskip 75pt=  \hat\psi(0) - \foh \psi(0), &  {\rm if} \; \rho(\phi) = \Ad\phi.
\end{cases}
\label{Wtransform}\end{equation}

\end{comment}

\subsection{Outline}\

The rest of the paper is organized as follows. In \S\ref{sec:preAut}, we collect various preliminaries on automorphic forms on $\GL_{3}(\Z)$ (their Fourier development and $L$-functions), and the Kontorovich-Lebedev-Whittaker transform, as explicated by the authors in \cite{GoldfeldKontorovich2011}. 
In \S\ref{sec:preKuz}, we collect the $\GL(3)$ Kuznetsov formula, explicating all the terms which appear. 

The careful definition of the choice of test function is given in \S\ref{sec:pIs}, where we also %carefully
 analyze its growth/decay properties; this is the most important and involved section. We note that, though the argument is a bit complicated (four-dimensional integrals of 12 Gamma factors in the numerator and 7 Gamma factors in the denominator),  the analysis uses nothing more than Stirling's asymptotics for the Gamma function.
In \S\ref{sec:boundsKloos}, we input the estimates of \S\ref{sec:pIs} into the Kloosterman integrals appearing on the geometric side of the Kuznetsov formula, giving bounds for these, as well as estimating away the contribution from the Eisenstein spectrum.
Combining all the above estimates, we prove Theorem \ref{thm:Kuz} in \S\ref{sec:PfThm1}. 

Next we turn our attention to the application to low-lying zeros. In \S\ref{subs:ExpForm}, we develop the Explicit Formula for the various $L$-functions of interest, and analyze the local Langlands-Satake parameters in \S\ref{subs:Loc}. Having done so, we apply Theorem \ref{thm:Kuz} to the low-lying zeros sum in \S\ref{sec:pfMain} to prove Theorem \ref{thm:main}.

\subsection*{Acknowledgements}\

The authors wish to thank Peter Sarnak for suggesting the application of our work on the Kontorovich-Lebedev transform to low-lying zeros. We are grateful to him and
 Valentin Blomer   for many discussions, comments and suggestions regarding this work. Thanks also to Steve J. Miller and Matt Young for comments on an earlier draft. Much of this work was carried out during the 2009-2010 special year in analytic number theory at IAS, and it is a pleasure for the authors to acknowledge the fantastic working conditions. This work has its roots dating back to the AIM workshop ``Analytic theory of GL(3) automorphic forms and applications'' in November 2008, and we also thank the organizers of this meeting.

%%%%%%%%%%%%%%%%%%%%%%%%%%%%%%%%%%%%%%%%
%%%%%%%%%%%%%%%%%%%%%%%%%%%%%%%%%%%%%%%%

\newpage

\section{Preliminaries on Automorphic Forms on $\GL_3(\Z)$}\label{sec:preAut}

\subsection{Jacquet's Whittaker Function}\label{subs:Whit}\
\vskip 5pt
%We follow \cite{Goldfeld2006}. Let
Let
$$
\frak h^3
:=
\GL_{3}(\R)/(O_{3}(\R)\times \R^{\times})
%\cong
%\SL_{3}(\R)/\SO_{3}(\R)
$$
denote  the generalized upper half plane. 
For $z\in\fh^{3}$ we use Iwasawa coordinates:
$$
z=xy=
\left(
\begin{array}{ccc}
1 & x_2 & x_3 \\
0 & 1 & x_1 \\
0 & 0 & 1\\
\end{array}\right) \left(
\begin{array}{ccc}
y_1y_2& 0 & 0 \\
0 & y_1 & 0 \\
0 & 0 & 1\\
\end{array}\right) 
,
$$
where $x_{1},x_{2},x_{3}\in\R$ and $y_{1},y_{2}>0$. We will frequently abuse notation, not distinguishing between $y$ as above and $y=(y_{1},y_{2})$.
Equip $\fh^{3}$ with the Haar measure
$$
dz={dx_{1}dx_{2}dx_{3}dy_{1}dy_{2}\over (y_{1}y_{2})^{3}}.
$$
With this measure, the group $\G=\GL_{3}(\Z)$ is a lattice, that is, the quotient $\G\bk\fh^{3}$ has finite volume. In fact, the volume is
\be\label{eq:volIs}
\int\limits_{\G\bk\fh^{3}}dz
=
{3
\gz(3)
\over
2\pi}.
\ee

For the pair $\nu=(\nu_{1},\nu_{2})\in\C^{2}$, set
\be\label{eq:nu3Def}
\nu_{3}:=\nu_{1}+\nu_{2}.
\ee
Then we have the $I$-function, defined by
\be\label{eq:Idef}
I_\nu(z) = (y_1y_{2})^{1+\nu_3}y_{1}^{\nu_2} y_2^{\nu_1}, \qquad (z = xy \in \frak h^3).
\ee
We now define Jacquet's Whittaker function for $GL_{3}(\mathbb R).$\footnote{Throughout we use the {\it completed} Whittaker function, in the terminology of \cite{Goldfeld2006}.}

\begin{definition}[Whittaker function]
%Let $U_3(\mathbb R)$ denote the unipotent group of upper triangular matrices in $GL_{3}(\mathbb R).$ Fix $\nu = \{\nu_1, \nu_2\} \in \Bbb C^2.$ 
For $\nu \in \C^2$ and $z \in \frak h^3$, set
% the  Whittaker function $\cW_\nu^{\pm}$ by the integral
\beann
W^{\pm}_\nu(z) 
&:=& 
 \pi^{-3\nu_{3}}\prod_{j=1}^{3}\Gamma\left(\frac{1+3\nu_j}{2}   \right) % \Gamma\left(\frac{3\nu_2}{2}   \right) \Gamma\left(\frac{3\nu_0-1}{2}   \right)
\\
&&
\times
\iiint\limits_{\R^{3}}
I_{\nu}\left(
\bp
&&1\\
&1&\\
1&&
\ep
\bp
1&u_{2}&u_{3}\\
&1&u_{1}\\
&&1
\ep
 z\right)
 e(-u_{1}\mp u_{2})
  \, du_{1}du_{2}du_{3}.
\eeann
\end{definition}
This function, originally defined for $\Re\nu_{1}, \Re\nu_{2} \gg1$, has analytic continuation to all $\nu\in\C^{2}$. 
For $z=y$, the value of $W^{\pm}_{\nu}(y)$ is independent of the sign, so we drop the $\pm$.

It is convenient to define the parameters $\ga$, given in terms of  $\nu$,  by the following linear relation:
\be\label{eq:ajDef}
\ga_{1}=2\nu_{1}+\nu_{2},\quad
\ga_{2}=-\nu_{1}+\nu_{2},\quad
\ga_{3}=-\nu_{1}-2\nu_{2}.
\ee
Then $\alpha_1 + \alpha_2 + \alpha_3 = 0$
and there is an action of the Weyl group which permutes the parameters $\alpha_1, \alpha_2, \alpha_3.$ We say that a  function of $\nu_1, \nu_2$ is symmetric under the action of the Weyl group
if it is invariant under all reorders of the triple $(\alpha_1, \alpha_2, \alpha_3).$

Consider the representation of the Whittaker function as a double inverse Mellin transform \cite{Stade2001}
\be
\label{eq:invMellin}
W_{\nu}(y)
=%&=&
{y_{1}y_{2}\pi^{3/2}%\prod_{j=1}^{3}\G\left({1+3\nu_{j}\over2}\right) 
\over (2\pi i)^{2}}
%\\
%\nonumber
%&&
%\times
\int\limits_{(C_{1})}\int\limits_{(C_{2})}
{
\prod\limits_{j=1}^{3}
\G\left({s_{1}+\ga_{j}\over2}\right)
\G\left({s_{2}-\ga_{j}\over2}\right)
\over
4\pi^{s_{1}+s_{2}}
\G\left({s_{1}+s_{2}\over 2}\right)
}
y_{1}^{-s_{1}}
y_{2}^{-s_{2}}
ds_{1}
ds_{2}
,
\ee
for any $C_{1},C_{2}>0$. Here we use the standard convention that for $C \in \R$, the symbol $(C)$ denotes the line $C + i \R.$ Note that $W_{\nu}$ is symmetric under the action of the Weyl group.

For $s\in\C$, Stade's formula \cite{Stade2002} gives
\be\label{eq:Stade}
\iint\limits_{\R_{+}^{2}}
W_{\nu}(y)
\overline{W_{\mu}(y)}
(\det y)^{s}
{dy_{1}dy_{2}\over (y_{1}y_{2})^{3}}
=
{
\pi^{3(1-s)}
\over
\G\left({3s\over 2}\right)
}
\prod_{1\le j,k\le3}
\G\left({s+\ga_{j}+\overline{\gb_{k}}\over2}\right)
,
\ee
where $\mu_{3}=\mu_{1}+\mu_{2}$ and $\gb_{1},\gb_{2},\gb_{3}$ are defined in terms of $\mu_{1},\mu_{2}$ as in \eqref{eq:ajDef}. The left side above is originally only defined for $\Re(s)$ sufficiently large; of course the right side gives its meromorphic continuation.
\

\subsection{Kontorovich-Lebedev transform}\label{subs:KLt}\
\vskip 5pt
%\begin{definition}[Kontorovich-Lebedev transform] 
Next, we give the analogue of the Kontorovich-Lebedev transform for $\GL(3)$, often referred to as the Lebedev-Whittaker transform \cite{GoldfeldKontorovich2011, Wallach1992}.
%Let
%$y = \text{diag}(y_1y_2, y_1, 1) \in GL(3, \mathbb R)$.
Let $f:\R_+^{2}\to\C$ and define $f^\sharp:\C^{2}\to \C$ by
\be\label{eq:LWTrans}
f^\sharp(\nu)
:=
\iint\limits_{\R_{+}^{2}} f(y) W_{\nu}(y)\, \frac{dy_1 dy_2}{(y_1y_2)^3}, %\qquad (t = \{t_1, t_2\} \in \mathbb R^2)
\ee
provided the integral converges absolutely. 
Then $f^\sharp$ is termed the Lebedev-Whittaker transform of $f$.
Note that $f^{\sharp}$ inherits the property that it is symmetric under the action of the Weyl group.%\end{definition}

The inverse  transform is given as follows.
Assuming $g$ is invariant under the action of the Weyl group and has sufficient decay% (see \cite{GoldfeldKontorovich2011})
, we define
%\begin{thm}[Kontorovich-Lebedev inversion]
%There is some constant $c_r>0$ (probably $c_r=\pi^{-r}$) such that
\be
\label{eq:LWInv}
g^{\flat}(y):={1\over (\pi i)^2} \int\limits_{-i\infty}^{i\infty}\int\limits_{-i\infty}^{i\infty}g(\nu)
\overline{W_{\nu}(y)}
{d\nu_1 d\nu_2
\over \prod\limits_{j=1}^3
 \Gamma\left(
{3\nu_j \over 2}
\right)  \Gamma\left(
{-3\nu_j \over 2}
\right)}\;
.
\ee
%where we have defined $t_0 := t_1+t_2.$
%\end{thm}

A sufficient condition 
on the test functions above
(see \cite{GoldfeldKontorovich2011})
 is that $g(\nu)$ have holomorphic extension to a strip $-\eta <|\Re(\nu_{1})|,|\Re(\nu_{2})|<\eta$ (for some $\eta > 0$) and
in this strip satisfy
\be\label{eq:gNeeds}
|g(\nu)|< \exp\left(-\frac{3\pi}4\sum_{k=1}^{3}|\nu_{k}|\right)\prod_{k=1}^{3}(1+|\nu_{k}|)^{-10}.
\ee

Then under these growth assumptions  we have
$$
g=f^{\sharp}
\qquad
\Longleftrightarrow
\qquad
f=g^{\flat}
,
$$
and the Parseval-type relation:
%\begin{thm}[K-L Parseval]
\be\label{eq:KLParseval}%gin{align*}
 \iint\limits_{\R_{+}^{2}} f_{1}(y)
 \overline{f_{2}(y)}\, \frac{dy_1 dy_2}{(y_1y_2)^3}
 = 
\frac1{(\pi i)^{2}}\int\limits_{-i\infty}^{i\infty}  \int\limits_{-i\infty}^{i\infty} 
f_{1}^\sharp(\nu)
\overline{f_{2}^\sharp(\nu) }
{d\nu_1 d\nu_2
\over 
 \prod\limits_{j=1}^3
 \Gamma\left(
{3\nu_j \over 2}
\right)  \Gamma\left(
{-3\nu_j \over 2}
\right)}
.
\ee%nd{align*}
%\end{thm}

\

\subsection{Cusp Forms}\label{subs:cusp}\
\vskip 5pt
Take a Hecke-normalized basis of Maass cusp forms $\{\phi_{j}\}_{j=1,2,3\ldots}$ for $\cH_{0}$, the cuspidal subspace of $L^{2}(\G\bk \fh^{3})$. 
The form $\phi_{j}$ is of type $(1/3+\nu_{1}^{(j)},1/3+\nu_{2}^{(j)})$ where $\nu^{(j)}=(\nu_{1}^{(j)},\nu_{2}^{(j)})\in\C^{2}$ denote the spectral parameters. 
When speaking of a fixed Maass form $\phi$, we drop the superscript $(j)$.
% Having assumed the Ramanujan conjectures, 
For tempered forms,
the spectral parameters $\nu_{1}$ and $\nu_{2}$ are purely imaginary.
Then 
with $\nu_{3}=\nu_{1}+\nu_{2}$,
the Laplace eigenvalue $\lambda_\phi$  is related to $\nu$ by
$$
\gl_{\phi}=1-6(\nu_{1}^{2}+\nu_{2}^{2}+\nu_{3}^{2})
.
$$

 Each such $\phi$ has the Fourier-Whittaker development given by %Piatetski-Shapiro
  \cite{Shalika1973, PS1975, Goldfeld2006}:
 \be\label{eq:phiFourier}
\phi(z) = \sum_{\g\in U_2(\Z)\bk \SL_{2}(\Z)} \sum_{k_1\ge1}\sum_{k_{2}\neq0} 
{A_\phi(k_1,k_2)
\over 
 k_1 |k_2|
} \,
W_{\nu}^{\sgn(k_{2})}\left(
\bp
k_{1}|k_{2}|&&\\
&k_{1}&\\
&&1
\ep
\bp
\g & \\
 & 1\\
\ep
z
\right)
,
 \ee
with the Hecke normalization $A_{\phi}(1,1)=1$.

The $L$-function attached to $\phi$ is given by
$$
L(s,\phi):=\sum_{n\ge1}{A(1,n)\over n^{s}}
,
$$
where $A(1,n) = A_\phi(1,n)$, i.e., we have dropped the $\phi$ from the notation.
This constitutes a degree $3$ $L$-function, which in completed form has Euler product
$$
\gL(s,\phi%_j
):= \prod_p L_p(s,\phi%_j
) 
$$
 with local factors for $p<\infty$ of type
$$
L_p(s,\phi%_j
):= \prod_{k=1}^3\left(1-{\ga_k(p)\over p^s}\right)^{-1} = \left(1-{A%_j
(p,1)\over p^s}+{%\overline
{A%_j
(1,p)}\over p^{2s}}-{1\over p^{3s}}\right)^{-1} 
$$
and for $p=\infty$,
$$
L_\infty(s,\phi%_j
):= \pi
^{-\frac{3s}{2}} \prod_{k=1}^3\G\left({s+\ga_k\over 2}\right)^{-1}.
$$

The Rankin-Selberg $L$-function is
$$
L(s,\phi\times\tilde\phi)
:=
\gz(3s)
\sum_{k_{1},k_{2}}
{|A(k_{1},k_{2})|^{2}
\over
(k_{1}^{2}k_{2})^{s}}
.
$$
This $L$-function has a pole at $s=1$.
Standard Rankin-Selberg theory, together with Stade's formula \eqref{eq:Stade} shows that the $j$-th Maass form $\phi_{j}$
has $L^{2}$ norm given by
\be\label{eq:phiNorm}
\|\phi_{j}\|^{2}
=
6\
\cL_{j}\cdot
\prod_{k=1}^{3}
\G\left({1+3\nu_{k}^{(j)}\over 2}\right)
\G\left({1-3\nu_{k}^{(j)}\over 2}\right)
,
\ee
where
\be\label{eq:LjDef}
\cL_{j}:=\Res_{s=1}L(s,\phi_{j}\times\tilde\phi_{j}).
\ee

\

%%%%%%%%%%%%%%%%%%%%%%%%%%%%%%%%%%%%%%%%
%%%%%%%%%%%%%%%%%%%%%%%%%%%%%%%%%%%%%%%%

\newpage

\section{The $\GL(3)$ Kuznetsov Formula}\label{sec:preKuz}

%%%%%%%%%%%%%%%%%%%%%%%%%%%%%%%%%%%%%%%%
%%%%%%%%%%%%%%%%%%%%%%%%%%%%%%%%%%%%%%%%

%\newpage

%\subsection{Overview}\label{subs:Kuz}\

The following equation  is %an overview of 
the $\GL(3)$ Kuznetsov formula, as compiled from \cite{BumpFriedbergGoldfeld1988} and \cite{Blomer2011}:
\be\label{eq:Kuz}
\boxed{\cC+\cE_{min}+\cE_{max}%+\cE_{res}
=
\cM
+
\cK
+
\widetilde{\cK}+
\widetilde\cK^{\scriptscriptstyle \vee}
,}
\ee
where each component is explicated below. Let  $p:\R_+^2\to\C$ be a  test function with suitable decay properties; 
a sufficient condition is that
\be\label{eq:pBndAt0}
|p(y_{1},y_{2})|
%,|q(y_{1},y_{2})|
\ll (y_{1}y_{2})^{2+\vep}
,
\ee 
as $y_{1},y_{2}\to0$, and that $p$ % and $q$ 
is otherwise bounded.
Fix positive integers $n_{1}, n_{2}, m_{1}, m_{2}.$

The
 left hand side of \eqref{eq:Kuz}, called the spectral side, consists of cuspidal and Eisenstein contributions. The cuspidal contribution is given by
\be\label{eq:CuspIs}
\cC
=%&=&
\sum_{j}
A_{j}(m_{1},m_{2})
\overline{A_{j}(n_{1},n_{2})}
{|p^{\sharp}(\nu^{(j)}_{1},\nu_{2}^{(j)})|^{2}
\over 6\ \cL_{j}\
\prod\limits_{k=1}^{3}
\G\left({1+3\nu_{k}^{(j)}\over 2}\right)
\G\left({1-3\nu_{k}^{(j)}\over 2}\right)
}
,
\ee%ann
where the sum on $j$ is over cuspidal Hecke-Maass forms $\phi_{j}$ on $\GL(3, \R)$.
The minimal Eisenstein series contributes
\beann
\cE_{min}
&=&
{1%c_{1}
\over (4\pi i)^{2}}
\int\limits_{-i\infty}^{i\infty}%\nu_{1}\in  i\R}
\int\limits_{-i\infty}^{i\infty}%\nu_{2}\in  i\R}
A_{\nu}(m_{1},m_{2})
\overline{A_{\nu}(n_{1},n_{2})}
{
|p^{\sharp}(\nu_{1},\nu_{2})|^{2}
\over \prod\limits_{k=1}^{3}
\left|
\gz(1+3\nu_{k})
\G\left({1+3\nu_{k}\over 2}\right)
\right|^{2}
%\prod_{k=1}^{3}
%
%\G\left({-3\nu_{k}^{(j)}\over 2}\right)
}\;
d\nu_{1}
d\nu_{2}
,
\eeann
where % $c_{1}$ is an absolute constant, and 
the minimal Eisenstein coefficients satisfy
\be\label{eq:AnuBnd}
|A_{\nu}(n_{1},n_{2})|
\ll_{\vep}
(n_{1}n_{2})^{\vep}
.
\ee

Lastly, the maximal Eisenstein contribution is
\beann
\cE_{max}
&=&
{c\over 2\pi i}
\sum_{j=1}^\infty\;
\int\limits_{-i\infty}^{i\infty}%\nu\in  i\R}
{B_{\nu,r_{j}}(m_{1},m_{2})
\overline{B_{\nu,r_{j}}(n_{1},n_{2})}
\over
L(1, \Ad u_{j})
|L(1+3\nu, u_{j})|^{2}
}\;
{
\left|p^{\sharp}_{T,R}\left(\nu- \frac{ir_{j}}{3}, \; \frac{2ir_{j}}{3}\right)\right|^{2}
\over 
\left|
\G\left({1+3\nu-ir_{j}\over2}\right)
\G\left({1+2ir_{j}\over2}\right)
\G\left({1+3\nu+ir_{j}\over2}\right)
\right|^{2}
}\;
d\nu,
\eeann
where $c$ is an  absolute constant, and $\{u_{j}\}$ is a basis of Hecke-Maass forms for $\GL(2,\Z)$, each of eigenvalue $1/4+r_{j}^{2}$.
The trivial bound for these Fourier coefficients is
\be\label{eq:BnuBnd}
|B_{\nu,r_{j}}(n_{1},n_{2})|
\ll_{\vep}
(n_{1}n_{2})^{1/2+\vep}
. %,
\ee
%see \cite[Proposition 10.9.3]{Goldfeld2006}.
%since we have assumed the Ramanujan conjectures (for $\GL(2)$).
%
%
Note that the residual spectrum does not contribute, having only degenerate terms in its Fourier expansion.

For  functions $p, q:\R_+^2\to\C$ let
$$\<p, q\> = \iint\limits_{\R_+^2} p(y_1, y_2) \,\overline{q(y_1, y_2)}\; \frac{dy_1 dy_2}{(y_1 y_2)^3}.$$

 Let $\bo_C$ denote the % Kronecker Delta 
 indicator function,
  which is $1$ if the % some %arithmetic 
 condition $C$ holds and %is 
 $0$ otherwise. 
The right-hand side of \eqref{eq:Kuz}, called the arithmetic side of the Kuznetsov  formula, consists of a main term and Kloosterman contributions given by
\bea
\label{eq:cMis}
\cM
&=&
\bo_{\left\{n_{1}=m_{1}\atop n_{2}=m_{2}\right\}}
\<p,p\>,
\\
\nonumber
\widetilde\cK
&=&
\sum_{\gep=\pm1}
\sum_{D_{1}\mid D_{2}\atop m_{2}D_{1}^{2}=n_{1}D_{2}}
{
\widetilde S(\gep m_{1},n_{1},n_{2},D_{1},D_{2})
\over
D_{1}D_{2}
}
\widetilde\cJ_{\gep}%_{p,q}
\left(\sqrt{n_{1}n_{2}m_{1}\over D_{1}D_{2}} %; \, \gep
\right),
\\
\nonumber
\widetilde\cK^{\scriptscriptstyle \vee}
&=&
\sum_{\gep=\pm1}
\sum_{D_{2}\mid D_{1}\atop m_{1}D_{2}^{2}=n_{2}D_{1}}
{
\widetilde S(\gep m_{2},n_{2},n_{1},D_{2},D_{1})
\over
D_{1}D_{2}
}
\widetilde\cJ_{\gep
%p, q
}%^{\scriptscriptstyle \vee}
\left(\sqrt{n_{1}n_{2}m_{2}\over D_{1}D_{2}} %;\, \gep
\right),
\\
\nonumber
\cK
&=&
\sum_{\gep_{1},\gep_{2}=\pm1}
\sum_{D_{1},D_{2}}
{
S(\gep_{1}m_{1},\gep_{2}m_{2},n_{1},n_{2},D_{1},D_{2})
\over
D_{1}D_{2}
}
\cJ_{\gep_{1},\gep_{2}%p, q
}
\left(
{
\sqrt{m_{1}n_{2}D_{1}}
\over
D_{2}
}
,
{
\sqrt{m_{2}n_{1}D_{2}}
\over
D_{1}
}
%;\, \gep_{1},\gep_{2}
\right)
.
\eea
Here $S, \widetilde S$, $\cJ, \widetilde\cJ$ are certain $\GL(3)$ Kloosterman sums and integrals corresponding to various elements of the Weyl group.

Let $e(x) := e^{2\pi i x}.$ The Kloosterman sums  are given explicitly by:
\beann
\widetilde S
(m_{1},n_{1},n_{2},D_{1},D_{2})
&:=&
\bo_{D_{1}\mid D_{2}}
\underset
{
{
C_{1}(\mod D_{1})
,
C_{2}(\mod D_{2})
}
\atop
{
(C_{1},D_{1})
=1
=(C_{2},D_{2}/D_{1})
}
}
{\sum\sum}
e\left(
{
m_{1}C_{1}+n_{1}\bar C_{1}C_{2}
\over
D_{1}
}\right)
e\left(
{
n_{2}\bar C_{2}
\over
D_{2}/D_{1}
}\right)
,
\eeann
%%%%%%%%%%%%%%%%%%%%%%%%%%%%
and
%%%%%%%%%%%%%%%%%%%%%%%%%%%%%
\beann
S(m_{1},m_{2},n_{1},n_{2},D_{1},D_{2})
&:=&
\underset
{
{
B_{1},C_{1}(\mod D_{1})
\atop
B_{2},C_{2}(\mod D_{2})
}
\atop
{
(B_{1},C_{1},D_{1})
=1
=(B_{2},C_{2},D_{2})
\atop
B_{1}B_{2}+C_{2} D_{2}+C_{2}D_{1}\equiv0(\mod D_{1}D_{2})
}
}
{\sum\sum\sum\sum}
e\left(
{
m_{1}B_{1}+n_{1}(Y_{1}D_{2}-Z_{1}B_{2})
\over
D_{1}
}\right)
\\
&&
\hskip1.7in
\times\
e\left(
{
m_{2}B_{2}+n_{2}(Y_{2}D_{1}-Z_{2}B_{1})
\over
D_{2}
}\right)
,
\eeann
where $Y_{1},Y_{2},Z_{1},Z_{2}$ are determined by
$$
Y_{1}B_{1}+Z_{1}C_{1}\equiv1 \, (\mod D_{1})
\qquad\text{and}\qquad
Y_{2}B_{2}+Z_{2}C_{2}\equiv1 \, (\mod D_{2})
.
$$

\

The Kloosterman integrals are given by:
\beann
\widetilde\cJ_{\gep%p,q
}
(A%; \gep
)
&=&
A^{-2}
\iint\limits_{\R_{+}^{2}}
\iint\limits_{\R^{2}}
\overline{
p(Ay_{1},y_{2})
}\,
e(-\gep A x_{1}y_{1})\
p\left(
{y_{2}}\cdot{
{\sqrt{1 +  x_{1}^{2}+x_{2}^{2} }}\over 1+ x_{1}^{2}}
, \;\;
\frac{A}{y_{1}y_{2}}\cdot{
\sqrt{1 + x_{1}^{2}}\over 1+ x_{1}^{2}+x_{2}^{2}}
\right)
\\
&&
\hskip1in
\times\
e\left(
{y_{2}}\cdot{
x_1x_2\over 1 + x_{1}^{2}} \;
+ \;
\frac{A}{y_{1}y_{2}}\cdot{
x_{2}\over 1 + x_{1}^{2}+x_{2}^{2}}
\right)
dx_{1}dx_{2}{dy_{1}dy_{2}\over y_{1}y_{2}^{2}}
,
\eeann
%%%%%%%%%%%%%%%%%%%%%%%%%%%%%%
%%%%%%%%%%%%%%%%%%%%%%%%%%%%%%
and
%%%%%%%%%%%%%%%%%%%%%%%%%%%%%%
\beann
\cJ_{\gep_{1},\gep_{2}
%p,q
}(A_{1},A_{2}%; \gep_{1},\gep_{2}
)
&=&
(A_{1}A_{2})^{-2}
\iint\limits_{\R_{+}^{2}}
\iiint\limits_{\R^{3}}
\overline{p(A_{1}y_{1},A_{2}y_{2})}
e(
-\gep_{1}A_{1}x_{1}y_{1}
-\gep_{2}A_{2}x_{2}y_{2}
)
\\
&&
\times\
p\left(
\frac{A_{2}}{y_{2}}\cdot{
\sqrt{(x_{1}x_{2}-x_{3})^{2}+x_{1}^{2}+1}\over x_{3}^{2}+x_{2}^{2}+1}
,
\frac{A_{1}}{y_{1}}\cdot{
\sqrt{x_{3}^{2}+x_{2}^{2}+1}\over (x_{1}x_{2}-x_{3})^{2}+x_{1}^{2}+1}
\right)
\\
&&
\times\
e\left(
-\frac{A_{2}}{y_{2}}\cdot{x_{1}x_{3}+x_{2}\over x_{3}^{2}+x_{2}^{2}+1}
-\frac{A_{1}}{y_{1}}\cdot{x_{2}(x_{1}x_{2}-x_{3})+x_{1}\over (x_{1}x_{2}-x_{3})^{2}+x_{1}^{2}+1}
\right)
\\
&&
\hskip3in
\times\
dx_{1}
dx_{2}
dx_{3}
{dy_{1}dy_{2}\over y_{1}y_{2}}
.
\eeann
%%%%%%%%%%%%%%%%%%%%%%%%%%%%%%%%%%%%%%%%
%%%%%%%%%%%%%%%%%%%%%%%%%%%%%%%%%%%%%%%%

\

\newpage

\section{Choice of Test Function and Bounds}\label{sec:pIs}\vskip 5pt

We now make a specific choice for the test function  $p(y_{1},y_{2})$.  
By Lebedev-Whittaker inversion \eqref{eq:LWInv}, 
we can just as well choose the transform $p^{\sharp}(\nu_{1},\nu_{2})$. Let $R \ge 10$ and  $T\gg 1.$ We define
\be\label{eq:pSharpIs}
p_{T,R}^{\sharp}(\nu_{1},\nu_{2})
:=
\sqrt6\
e^{\ga_{1}^{2}+\ga_{2}^{2}+\ga_3^2\over 2T^{2}}
\prod_{1\le j\le 3} 
\G\left({2+R+3\nu_j\over4}\right)
\G\left({2+R-3\nu_j\over4}\right)
.
\ee
This choice is motivated by the fact that we need $p^{\sharp} = p^{\sharp}_{T,R}$ to be invariant under the action of the Weyl group, while also  
requiring 
 cancellation of  the exponential growth of the $\G$-factors in the denominator on the right side of \eqref{eq:CuspIs} (cuspidal contribution to the Kuznetsov  formula). The variable $R\ge 10$ is introduced to obtain absolute convergence of the sum \eqref{eq:CuspIs}, and to pull certain contours without passing through poles, see  \eqref{eq:RtoGk}.
Note first that $p_{T,R}^{\sharp}$ easily satisfies the requisite bounds \eqref{eq:gNeeds} for Lebedev-Whittaker inversion.
It will be shown below that the inverse transform $p_{T,R}$ satisfies \eqref{eq:pBndAt0}, see \eqref{eq:range2}.

\vskip 5pt
Observe then that the cuspidal contribution  \eqref{eq:CuspIs} becomes
\be
\cC=
\sum_{j}
{A_{j}(m_{1},m_{2})
\overline{A_{j}(n_{1},n_{2})}
}
{h_{T,R}(\nu_{1}^{(j)},\nu_{2}^{(j)})\over \cL_{j}}
,\notag
\ee
exactly as desired in \eqref{eq:hTdef}.
\vskip 10pt

%\newpage

\subsection{Some Auxiliary Bounds}\
\vskip 5pt

 We collect here some bounds coming from Stirling's asymptotic formula:
 %
 %we shall need Stirling's asymptotic formula for the Gamma function 
%   which takes the form 
\begin{equation}  \label{Stirling} |\Gamma(\sigma+it)| \sim \sqrt{2\pi}\, |t|^{\sigma-\frac12}\, e^{-\frac{\pi|t|}{2}}, \qquad (t \to \pm\infty),\end{equation}for fixed $\sigma \in \R$.

%  (\ref{Stirling}) for the absolute value of the Gamma function.
\vskip 10pt
There are three types of integrals which we will need to estimate $p_{_{T,R}}$. 
%As always, $\ga$'s are related to $\nu$'s by 
%\eqref{eq:gaIs}. 
Throughout we have $y_{1},y_{2}>0,$ $R\ge10$, $T\gg 1$. 
%We adopt the notation that $\int_{(c)} := \int_{c-i\infty}^{c+i\infty}$ for $c\in \R.$

\vskip 10pt

\noindent
$\underline{\text{\bf The First Integral:}}$ 
\vskip 5pt
For any $C_{1},C_{2}\in\R\setminus\{-2,-4,-6,\dots\}$, let
\bea
\label{eq:I1def}
\hskip-.5in
\cI_{T,R}^{(1)}(C_{1},C_{2};y_{1},y_{2})
&:=&
\int\limits_{(0)}
\int\limits_{(0)}
\int\limits_{(C_{2})}
\int\limits_{(C_{1})}
{
e^{\ga_{1}^{2}+\ga_{2}^{2}+\ga_3^2\over 2T^{2}}
\prod\limits_{1\le j\le 3}
\G\left(
{2+R+3\nu_j\over 4}
\right)
\G\left(
{2+R-3\nu_j\over 4}
\right)
\over
\prod\limits_{j=1}^{3}
\G
\left(
{3\nu_{j}\over 2}
\right)
\G
\left(
{-3\nu_{j}\over 2}
\right)
}
\\
\nonumber
&&
\hskip60pt
\times
{
\prod\limits_{j=1}^{3}
\G\left({s_{1}-\ga_{j}\over 2}\right)
\G\left({s_{2}+\ga_{j}\over 2}\right)
\over
4\pi^{s_{1}+s_{2}}
\G\left({s_{1}+s_{2}\over 2}\right)
}
y_{1}^{1-s_{1}}
y_{2}^{1-s_{2}}
ds_{1}
ds_{2}
d\nu_{1}
d\nu_{2}
.
\eea

\noindent
$\underline{\text{\bf The Second Integral:}}$
\vskip 5pt Similarly, for any $\gk_{1},\gk_{2},C_{1}\in\R$ (so that the integrand below doesn't pass through poles of $\G$), let
\beann
\hskip-.5in
\cI_{T,R}^{(2)}(\gk_{1},\gk_{2},C_{1};y_{1},y_{2})
&:=& \hskip -7pt
\int\limits_{(\gk_{2})}
\int\limits_{(\gk_{1})}
\int\limits_{(C_{1})}
{
e^{\ga_{1}^{2}+\ga_{2}^{2}+\ga_3^2\over 2T^{2}}
\prod\limits_{1\le j\le 3}
\G\left(
{2+R+3\nu_j\over 4}
\right)
\G\left(
{2+R-3\nu_j\over 4}
\right)
\over
\prod\limits_{j=1}^{3}
\G\left({3\nu_{1}\over 2}\right)
\G\left({3\nu_{2}\over 2}\right)
\G\left({3\nu_{3}\over 2}\right)
\G\left({-3\nu_{2}\over 2}\right)
}
\\
&&
\hskip35pt
\times
{
\G\left({s_{1}-\ga_{2}\over 2}\right)
\G\left({s_{1}-\ga_{3}\over 2}\right)
}
y_{1}^{1-s_{1}}
y_{2}^{1+\ga_{1}}
ds_{1}
d\nu_{1}
d\nu_{2}
.
\eeann

\noindent
$\underline{\text{\bf The Third Integral:}}$
\vskip 5pt
Lastly, for any $\gk_{1},\gk_{2}\in\R$ passing through no poles, let
\beann
\hskip-.5in
\cI_{T,R}^{(3)}(\gk_{1},\gk_{2};y_{1},y_{2})
&:=&
\int\limits_{(\gk_{2})}
\int\limits_{(\gk_{1})}
{
e^{\ga_{1}^{2}+\ga_{2}^{2}+\ga_3^2\over 2T^{2}}
\prod\limits_{1\le j\le 3}
\G\left(
{2+R+3\nu_j\over 4}
\right)
\G\left(
{2+R-3\nu_j\over 4}
\right)
\over
\G\left({3\nu_{1}\over 2}\right)
\G\left({3\nu_{3}\over 2}\right)
\G\left({-3\nu_{2}\over 2}\right)
}
\\
&&
\hskip176pt
\times
y_{1}^{1-\ga_{2}}
y_{2}^{1+\ga_{1}}
d\nu_{1}
d\nu_{2}.
\eeann

Define $\gk_{1}',\gk_{2}',\gk_{3}'$ to be related to $\gk$'s in the same way that $\ga$'s are related to $\nu$'s, that is, 
\be\label{eq:gkP}
\gk'_1=2\gk_1+\gk_2,
\quad
\gk'_2=-\gk_1+\gk_2,
\quad
\gk'_3=-\gk_1-2\gk_2.
\ee

\begin{thm} \label{thm:bnds} 
Fix $R\ge 10$ and $\vep > 0.$
For any $y_{1},y_{2}>0$ and $T\gg 1$, 
%\begin{enumerate}
%\item
we have the bound
\be\label{eq:I1}
|\cI_{T,R}^{(1)}(C_{1},C_{2};y_{1},y_{2})| \,
\ll_{\vep,C_{1},C_{2},R} \;
(y_1y_2)
 T^{9/2+3R/2}
\left(
y_{1}\over T
\right)^{-C_{1}}
\left(
y_{2}\over T
\right)^{-C_{2}}
T^{\vep}
.
\ee

%\item
Moreover,
\be\label{eq:I2}
|\cI_{T,R}^{(2)}(\gk_{1},\gk_{2},C_{1};y_{1},y_{2})| \,
\ll_{\vep,C_{1},\gk_{1},\gk_{2},R} \;
(y_1 y_2)
T^{ 4+3R/2} 
\left(
{
y_{1}\over T
}
\right)^{-C_{1}}
\left(
{y_{2}
\over
T}
\right)^{\gk'_{1}}
T^{\vep}
.
\ee

%\item
And finally,
\be\label{eq:I3}
|\cI_{T,R}^{(3)}(\gk_{1},\gk_{2};y_{1},y_{2})| \,
\ll_{\vep,\gk_{1},\gk_{2},R} \;
(y_1
y_2)
T^{7/2+3R/2}
\left(
{y_{1}\over T}
\right)^{-\gk'_{2}}
\left({y_{2}\over T}
\right)^{\gk'_{1}}
T^{\vep}
.
\ee

%\end{enumerate}
\end{thm}

We give separate treatments of each statement.

\pf[Proof of \eqref{eq:I1}]
\
 
Write $\nu_{j}=it_{j}$ and $s_{j}=C_{j}+iu_{j}$. The first exponential in the integrand gives arbitrary decay once $|t_{j}|>T^{1+\vep}$ for any $\vep$.
Bringing the absolute values inside and applying Stirling's asymptotic formula gives

\beann
\hskip-.5in
|\cI_{T,R}^{(1)}(C_{1},C_{2};y_{1},y_{2})|
&\ll_{C_{1},C_{2},R,\vep}&
y_1^{1-C_1} y_2^{1-C_2}
\iint\limits_{|t_{1}|,| t_{2}|\le T^{1+\vep}}
\iint\limits_{\R^{2}}
\cP
\cdot
\exp
\left(
\frac\pi4\cdot \cE
\right)
du_{1}
du_{2}
dt_{1}
dt_{2}
,
\eeann
where $\cE=\cE(t_{1},t_{2},u_{1},u_{2})$ is the exponential factor:
\beann
\cE
&=&
3\sum_{k=1}^{3}|t_{k}|
%+3 \left|t_1\right|+3 \left|t_2\right|+3 \left|t_1+t_2\right|+
-
\sum_{k=1}^{3}
|\ga_{k}-iu_{1}|
%-\left|2 t_1+t_2-u_1\right|-\left|t_1-t_2+u_1\right|-\left|t_1+2 t_2+u_1\right|
%\\
%&&
-
\sum_{k=1}^{3}
|\ga_{k}+iu_{2}|
%-\left|t_1+2 t_2-u_2\right|-\left|-t_1+t_2+u_2\right|-\left|2 t_1+t_2+u_2\right|
%\\
%&&
+\left|u_1+u_2\right|
,
\eeann
and  $\cP=\cP_{C_{1},C_{2},R}(t_{1},t_{2},u_{1},u_{2})$ is the polynomial factor:
\beann
\cP
&=&
%\left(\left(\left|t_1\right|+1\right) \left(\left|t_2\right|+1\right) \left(\left|t_1+t_2\right|+1\right)\right){}^{\frac{R+2}{2}} 
\left(\prod\limits_{k=1}^{3}(1+|t_{k}|)\right)^{(R+2)/2}
%\\
%&&
\left(\prod_{k=1}^{3}(1+|\ga_{k}-iu_{1}|)\right)^{(C_{1}-1)/2}
\\
&&
\times
%\left(\left(\left|2 t_1+t_2-u_1\right|+1\right) \left(\left|t_1-t_2+u_1\right|+1\right) \left(\left|t_1+2 t_2+u_1\right|+1\right)\right){}^{\frac{1}{2} \left(C_1-1\right)} 
\left(\prod_{k=1}^{3}(1+|\ga_{k}+iu_{2}|)\right)^{(C_{2}-1)/2}
%\\
%&&
%\times
%\left(\left(\left|t_1+2 t_2-u_2\right|+1\right) \left(\left|-t_1+t_2+u_2\right|+1\right) \left(\left|2 t_1+t_2+u_2\right|+1\right)\right){}^{\frac{1}{2} \left(C_2-1\right)}
%\\
%&&
%\times
\Bigg(1+\left|u_1+u_2\right|\Bigg)^{ \left(1-C_1-C_2\right)/2} 
.
\eeann

Note that we always have 
$$
\cE\le
0
,
$$ 
with equality only when
$$
t_2-t_1\leq u_1\leq 2 t_1+t_2
\qquad\text{and}\qquad
 t_1-t_2\leq u_2\leq t_1+2 t_2
$$
or
$$%\be\label{eq:ineq2}
-t_1-2 t_2\leq u_1<t_2-t_1
\qquad\text{and}\qquad
-2 t_1-t_2\leq u_2\leq t_1-t_2
 .
$$%\ee
Hence, there is %must be 
arbitrary decay outside of this range.
Both inequalities have the same contribution, so we only deal with the second.

Make a linear change variables 
$$
u_1\mapsto u_{1} -t_1-2 t_2
\qquad\text{and}\qquad
u_2\mapsto u_{2}-2 t_1-t_2
,
$$
so the new range is
\be\label{eq:ut}
0\leq u_1< 3 t_2
\qquad\text{and}\qquad
0\leq u_2\leq 3 t_1
,
\ee
and the $\cP$ factor becomes
\beann
\cP_{1}
&:=&
\left(\left(1+|t_1|\right) \left(1+|t_2|\right) \left(1+|t_1+t_2|\right)\right)^{(R+2)/2} 
\\
&&
\times
\left(\left(1+\left|u_1\right|\right) \left(1+\left|3 t_1+3 t_2-u_1\right|\right) \left(1+\left|u_1-3 t_2\right|\right)\right)^{ \left(C_1-1\right)/2} 
\\
&&
\times
\left(\left(1+\left|u_2\right|\right) \left(1+\left|3 t_1+3 t_2-u_2\right|\right) \left(1+\left|u_2-3 t_1\right|\right)\right)^{ \left(C_2-1\right)/2} 
\\
&&
\times
\left(1+\left|-3 t_1-3 t_2+u_1+u_2\right|\right)^{\left(1-C_1-C_2\right)/2}
.
\eeann

The integral of $\cP_{1}$ over \eqref{eq:ut} in $u_{1},u_{2}$ is bounded up to constant by
\beann
\cP_{2}
&:=&
\left(1+|t_1|\right)^{(R+2)/2+C_{2}} 
\left(1+|t_2|\right)^{(R+2)/2+C_{1}} 
\left(1+|t_1+t_2|\right)^{(R+1)/2} 
.
\eeann
%One can readily check that the second inequality \eqref{eq:ineq2} gives the same contribution.
%
%
Integrating $\cP_{2}$ over the range $|t_{j}|<T^{1+\vep}$ gives \eqref{eq:I1}, as claimed.
\epf

\

\

Next we give a 
\pf[Proof of \eqref{eq:I2}]
\

Again by Stirling's formula, we have
\beann
\hskip-.5in
|\cI_{T,R}^{(2)}(\gk_{1},\gk_{2},C_{1};y_{1},y_{2})|
&\ll_{C_{1},\gk_{1},\gk_{2},R,\vep}&
y_1^{1-C_1} y_2^{1+\gk'_{1}}
\iint\limits_{| t_{1}|,| t_{2}|\le T^{1+\vep}}
\int\limits_{\R}
\cP
\cdot
\exp
\left(
\frac\pi4\cdot \cE
\right)
du_{1}
dt_{1}
dt_{2}
,
\eeann
where $\cE=\cE(t_{1},t_{2},u_{1}%,u_{2}
)$ is now the exponential factor:
\beann
\cE
&=&
-\left|t_1-t_2+u_1\right|-\left|t_1+2 t_2+u_1\right|+3| t_2|
,
\eeann
and  $\cP=\cP_{C_{1},C_{2},\gk_{1},\gk_{2},R}(t_{1},t_{2},u_{1}%,u_{2}
)$ is now the polynomial factor:
\beann
\cP
&=&
\left(1+|t_2|\right){}^{\frac{R}{2}+1} 
\left(1+|t_1|\right){}^{\frac{1}{2} \left(-3 \kappa _1+R+1\right)} 
\left(1+|t_1+t_2|\right){}^{\frac{1}{2} \left(-3 \kappa _1-3 \kappa _2+R+1\right)} 
\\
&&
\times
\left(1+\left|t_1-t_2+u_1\right|\right){}^{\frac{1}{2} \left(C_1+\kappa _1-\kappa _2-1\right)} 
\left(1+\left|t_1+2 t_2+u_1\right|\right){}^{\frac{1}{2} \left(C_1+\kappa _1+2 \kappa _2-1\right)}
.
\eeann

Note that we always have 
$$
\cE\le
0
,
$$ 
with equality only when
$$
 -t_1-2 t_2\leq u_1\leq t_2-t_1
,
$$
so we may restrict the $u_{1}$ integral to this range.

Make a linear change variables 
$$
u_1\mapsto u_{1} -t_1-2 t_2
,
$$
so the new range is
$$
0\leq u_1< 3 t_2
,
$$
and the $\cP$ factor becomes
\beann
\cP_{1}
&:=&
\left(1+|t_2|\right){}^{(R+2)/2} 
\left(1+|t_1|\right){}^{ \left(-3 \kappa _1+R+1\right)/2} 
\left(1+|t_1+t_2|\right){}^{\left(-3 \kappa _1-3 \kappa _2+R+1\right)/2} 
\\
&&
\times
\left(\left|u_1\right|+1\right){}^{ \left(C_1+\kappa _1+2 \kappa _2-1\right)/2} \left(\left|u_1-3 t_2\right|+1\right){}^{\left(C_1+\kappa _1-\kappa _2-1\right)/2}
.
\eeann

The integral of $\cP_{1}$ over the $u_{1}$ range is bounded up to constant by
\beann
\cP_{2}
&=&
\left(1+|t_1|\right){}^{ \left(-3 \kappa _1+R+1\right)/2} 
\left(1+|t_2|\right){}^{\left(R+2C_1+2\kappa _1+\kappa _2+2\right)/2}
\\
&&
\times
\left(1+|t_1+t_2|\right){}^{\left(-3 \kappa _1-3 \kappa _2+R+1\right)/2} 
.
\eeann

Integrating $\cP_{2}$ 
over $|t_{j}|<T^{1+\vep}$
 gives
the claim.
\epf

Finally, we give a
\pf[Proof of \eqref{eq:I3}]
\

As before, we have
\beann
\hskip-.5in
|\cI_{T,R}^{(3)}(\gk_{1},\gk_{2};y_{1},y_{2})|
&\ll_{%C_{1},C_{2},
\gk_{1},\gk_{2},R,\vep}&
y_1^{1-\gk'_{2}} y_2^{1+\gk'_{1}}
\iint\limits_{%0\le 
|t_{1}|,| t_{2}|\le T^{1+\vep}}
\cP\
dt_{1}
dt_{2}
,
\eeann
where
$\cP=\cP_{%C_{1},C_{2},
\gk_{1},\gk_{2},R}(t_{1},t_{2}%,u_{1},u_{2}
)$ is the polynomial factor 
\beann
\cP
&=&
\left(1+|t_1|\right){}^{\left(-3 \kappa _1+R+1\right)/2} \left(1+|t_2|\right){}^{\left(3 \kappa _2+R+1\right)/2} \left(1+|t_1+t_2|\right){}^{\left(-3 \kappa _1-3 \kappa _2+R+1\right)/2}.
\eeann
(Note  that the exponential terms exactly cancel.)
Integrating $\cP$ gives the claim.
\epf

%\newpage

\subsection{Estimating $p_{_{T,R}}$}\label{subsec:pEst}\

\vskip 5pt
We use the bounds of the previous section to give an estimate for $p_{_{T,R}}$.
% in various ranges.
Among other things, we must verify that the inverse Lebedev-Whittaker transform $p_{T,R}$ satisfies \eqref{eq:pBndAt0}. This will follow from the bound \eqref{eq:range2}.

By Lebedev inversion \eqref{eq:LWInv}, we define
\be\label{eq:pTy}
p_{_{T,R}}(y)
:=
\frac1{(\pi i)^{2}}
\int\limits_{(0)}
\int\limits_{(0)} p^{\sharp}_{T,R}(\nu) \,\overline{W}_{\nu}(y){d\nu\over 
\prod\limits_{j=1}^{3}
\G
\left(
{3\nu_{j}\over 2}
\right)
\G
\left(
{-3\nu_{j}\over 2}
\right)
}
.
\ee
%The form of $q^{\sharp}_{T,R}$ means that the above integral is essentially cut off at $[-T,T]$.
%By playing with symmetries of the Weyl group, we may assume further (see Blomer (2.17)) that 
%\be\label{eq:t1t2Range}
%0\le t_{1}\le t_{2}\le T.
%\ee
%We need
%Also r

Recall
 the double inverse Mellin transform formula for the Whittaker function \eqref{eq:invMellin}, and that
 $\overline{W_{\nu}(y)}=W_{-\nu}(y)$ for $\nu$ tempered.

Then putting \eqref{eq:invMellin}  into \eqref{eq:pTy}
and comparing with 
\eqref{eq:I1def}
 gives
\be
\label{eq:pTRis}
\hskip-.5in
p_{_{T,R}}(y)
=%&=&
{\sqrt {6\pi^{3}}\over (\pi i)^{2}}\cdot
\cI_{T,R}^{(1)}(C_{1},C_{2};y_{1},y_{2})
,
\ee
and the equality holds  for any $C_{1},C_{2}>0$.
An immediate application of \eqref{eq:I1} proves that for any $y_{1},y_{2},C_{1},C_{2}>0$,
\be\label{eq:range0}
p_{_{T,R}}(y_{1},y_{2})
\ll_{C_{1},C_{2},
R,\vep}
%\boxed{
%\boxed{
y_1y_2 \,
 T^{9/2+3R/2}
\left(
y_{1}\over T
\right)^{-C_{1}}
\left(
y_{2}\over T
\right)^{-C_{2}}
T^{\vep}
%}}
.\ee

\begin{comment}
This gives us our first bound, using \eqref{eq:I1}.

\begin{lem}\label{lem:range0} Fix $R \ge 10$ and $\epsilon > 0.$
For any $C_{1},C_{2}>0$, $T\gg 1$, and
 any $y_{1},y_{2}>0$, we have
\end{lem}

This bound is optimal when $y_{1},y_{2}>T^{1+\vep}$.\\

\end{comment}

%Next % we return to \eqref{eq:pTRis} and 
\subsubsection{Pull past one set of poles}\

The above bound is insufficient for our purposes, so we return to
 the 
definition of $\cI^{(1)}_{T,R}$,
and
pull the $s_{2}$ integral from the vertical line $(C_{2})$ with $C_{2}>0$ to the vertical line $(-\cC_{2})$, with $\cC_{2}=-C_{2}$, $0<\cC_{2}<2$.
In so doing, we 
 pass through simple poles at $s_{2}=-\ga_{1},-\ga_{2},-\ga_{3}$ (generically the $\ga_{j}$ are distinct). Then we can write
 \be\label{eq:pOneSet}
 p_{_{T,R}}=\fM+\cR_{1}+\cR_{2}+\cR_{3}
 ,
 \ee
where $\fM$ is the remaining 4-dimensional integral
(that is, a constant times 
$\cI_{T,R}^{(1)}(C_{1},-\cC_{2};y_{1},y_{2})$),
 and the $\cR_{j}$ are the 3-dimensional contributions from the residues at $s_{2}=-\ga_{j}$. Note that $\cR_{1}$
%(the other two residues have similar contributions)
 is exactly equal to
 a constant times 
$\cI_{T,R}^{(2)}(0,0,C_{1};y_{1},y_{2})$. In this integral, 
we
 pull the $\nu_{1},\nu_{2}$ integrals from the vertical lines with $\Re\nu_{j}=0$ to the lines $(\gk_{1}),(\gk_{2})$ respectively,
 so that it becomes
$\cI_{T,R}^{(2)}(\gk_{1},\gk_{2},C_{1};y_{1},y_{2})$.
To ensure that we haven't passed any new poles, we require that the $\gk$'s satisfy:
\bea\label{eqs:gkToCs1}
|\kappa _j|&<&\frac{R+2}{3}
\\
\nonumber
\gk'_{2},\gk'_{3}&<&C_1
.
\eea
Recall that here, as always, the $\gk'$ are related to $\gk$ by \eqref{eq:gkP}. 

The estimate \eqref{eq:I2} %essentially 
bounds $\cR_{1}$ by 
$$
y_1 y_2\,
T^{ 4+3R/2} 
\left(
{
y_{1}\over T
}
\right)^{-C_{1}}
\left(
{y_{2}
\over
T}
\right)^{\gk'_{1}}
T^{\vep}
,
$$
whereas the term $\fM$ is dominated using \eqref{eq:I1} %essentially 
by
$$
y_1y_2\,
 T^{9/2+3R/2}
\left(
y_{1}\over T
\right)^{-C_{1}}
\left(
y_{2}\over T
\right)^{\cC_{2}}
T^{\vep}
.
$$
To 
%We can 
make these the same in $y_{2}$,
we would like to take $\gk_{1}'=2\gk_{1}+\gk_{2}$ as large as $\cC_{2}$, subject to \eqref{eqs:gkToCs1}, which requires $-\gk_{1}+\gk_{2}<C_{1}$ and $-\gk_{1}-2\gk_{2}<C_{2}$.
This is easily achieved by, say, setting $\gk_{2}=0, \gk_{1}>0$; then we can take $\gk_{1}$ as large as $1$, so that $\gk_{1}'$ can be as large as $2$.
We can take $\gk_{1}$ as large as $1$, as needed, since $R=10$.
%\newpage
%
% by choosing $\gk_{1}=\gk_{2}=\cC_{2}/3$, so that $\gk_{1}'=\cC_{2}$. The condition \eqref{eqs:gkToCs1} becomes
%$$
%\cC_{2}<R+2,\qquad
%\cC_{2}<C_{1}
%,
%$$
%which is satisfied (since $\cC_{2}<2$) if $R>0,$% and $C_{1}>2$
%. 
So under these conditions, we have dominated $\cR_{1}$ by the bound we already have on $\fM$. The same can be done with $\cR_{2}$ and $\cR_{3}$, by pulling $\gk$'s to different ranges.

We have thus given our second intermediate bound:
%Now we 
%So bounding $\fM$ using \eqref{eq:I1}, we summarize the above discussion with the following
%\begin{lem}
%Assume $R>1$.
%Fix $R\ge 10$ and $\vep > 0.$ For
for
any $y_{1},y_{2},C_{1}>0$, %$T\gg 1$, 
and
any
$0<\cC_{2}<2$,
%and 
% any $C_{1}>0%\cC_{2}
% $,
%any $\gk_{1},\gk_{2}\in\R$ satisfying \eqref{eqs:gkToCs1}, 
we have
\be\label{eq:range1}%$$%\beann
p_{_{T,R}}(y_{1},y_{2})
%&
\ll_{C_{1},\cC_{2},R,\vep}
%&
%\boxed{
%\boxed{
 y_1y_2\,
 T^{9/2+3R/2}
\left(
y_{1}\over T
\right)^{-C_{1}}
\left(
y_{2}\over T
\right)^{\cC_{2}}
T^{\vep}.
%}}
\ee%ann
%\end{lem}
%
By symmetry, we have the same result with the subscripts ``1'' and ``2'' reversed.

\begin{comment} 
By symmetry, we also have 
\begin{lem} Fix $R\ge 10$ and $\epsilon > 0.$
For
any $y_{1},y_{2}>0$, $T\gg 1$,
any
$0<\cC_{1}<2$,
and 
 any $C_{2}>0%\cC_{1}
 $,
%any $\gk_{1},\gk_{2}\in\R$ satisfying \eqref{eqs:gkToCs1}, 
%we have
\be\label{eq:range1p}%ann
p_{_{T,R}}(y_{1},y_{2})
%&
\ll_{\cC_{1},C_{2},R,\vep}
%&
\boxed{
\boxed{
y_1y_2\,
 T^{9/2+3R/2}
\left(
y_{1}\over T
\right)^{\cC_{1}}
\left(
y_{2}\over T
\right)^{-C_{2}}
T^{\vep}.
}}
\ee%ann
\end{lem}

These bounds apply optimally when one $y_{j}>T^{1+\vep}$ and the other is smaller.
\\

\end{comment}

\subsubsection{Pull past two sets of poles}\

The above is still insufficient, so  we return  to \eqref{eq:pOneSet}. 
In the $\fM$ integral, we now also pull the $s_{1}$ integral from $(C_{1})$ with $C_{1}>0$ to $(-\cC_{1})$, where $\cC_{1}=-C_{1}$, with $0<\cC_{1}<2$, passing through poles at $s_{1}=\ga_{j}$, giving
 $$
\fM=\tilde\fM+
 \tilde\cR_{1}+
 \tilde\cR_{2}+
 \tilde\cR_{3}
 .
 $$
The integral $\tilde \fM$ is exactly equal to an absolute constant times
$
\cI_{T,R}^{(1)}(-\cC_{1},-\cC_{2};y_{1},y_{2})
,
$
and
hence
we apply
 \eqref{eq:I1}, %essentially 
%bounds it by
giving
$$
\tilde \fM\ll_{\cC_{1},\cC_{2},R,\vep}
y_1y_2\,
 T^{9/2+3R/2}
\left(
y_{1}\over T
\right)^{\cC_{1}}
\left(
y_{2}\over T
\right)^{\cC_{2}}
T^{\vep}
.
$$
The integrals $\tilde\cR_{j}$ are of the same form as 
$
\cI_{T,R}^{(2)}(0,0,\cC_{2};y_{1},y_{2})
,
$
and after pulling to appropriate $\gk$'s, we can dominate the $\tilde\cR_{j}$ integrals by $\tilde\fM$, exactly as before.
\\

In the integral $\cR_{1}$, 
which is exactly
equal to a constant times
$
\cI_{T,R}^{(2)}(0,0,C_{1};y_{1},y_{2})
,
$
we can pull the $s_{1}$ integral from $(C_{1})$ to $(-\cC_{1})$, passing through poles at $s_{1}=\ga_{2},\ga_{3}$.
Hence we can write correspondingly
$$
\cR_{1}=\cR_{1}'+\cP_{1,2}+\cP_{1,3}.
$$
Here 
$\cR_{1}'$ is a triple integral, exactly equal to a constant times 
$
\cI_{T,R}^{(2)}(0,0,-\cC_{1};y_{1},y_{2})
,
$
and
$\cP_{1,2}$ is a double integral, exactly equal to a constant times
$
\cI_{T,R}^{(3)}(0,0;y_{1},y_{2})
.
$
The term $\cP_{1,3}$ is similar to $\cP_{1,2}$.

In the double integral $\cP_{1,2}$, we can pull contours in $\nu_{j}$ to any $(\gk_{j})$ with 
\be\label{eq:RtoGk}
|\gk_{j}|<(2+R)/3,
\ee
without passing new poles.
Apply the estimate \eqref{eq:I3} %to %essentially 
to bound $\cP_{1,2}$ by
$$
y_1
y_2\,
T^{7/2+3R/2}
\left(
{y_{1}\over T}
\right)^{-\gk'_{2}}
\left({y_{2}\over T}
\right)^{\gk'_{1}}
T^{\vep}
.
$$

%*********************** FIX THIS *************** gk1' $<$-$>$ gk2'

Elementary linear algebra shows from \eqref{eq:gkP} that if we choose 
$$
\gk_{1}=-\frac13(\cC_{1}+\cC_{2}),\qquad
\gk_{2}=\frac13(-\cC_{1}+2\cC_{2})
,
$$ 
then $-\gk_{1}'=\cC_{1}$ and $\gk_{2}'=\cC_{2}$. Since $0<\cC_{1},\cC_{2}<2$, the condition $R\ge10$ is more than sufficient to ensure that \eqref{eq:RtoGk} is satisfied.
So the $\cP_{1,2}$ contribution is dominated by that from $\tilde\fM$. The same (by a different pull in $\gk$'s) holds for $\cP_{1,3}$.

Lastly, consider the triple integral $\cR_{1}'
=
\cI_{T,R}^{(2)}(0,0,-\cC_{1};y_{1},y_{2})
.
$
 Since $-\cC_{1}=C_{1}<0$,  
 the next poles in $\nu$ arise when $\gk'_{2},\gk_{3}'=2+C_{1}$.
 Hence
 we can pull the $\nu$ variables to any $\gk_{1},\gk_{2}$, satisfying
\bea\label{eqs:gkToCs2}
|\kappa _j|&<&\frac{R+2}{3}
\\
\nonumber
\gk'_{2},\gk'_{3}&<&2-\cC_1
,
\eea
without passing more poles. 
The estimate \eqref{eq:I2} bounds $\cR'_{1}$ by
$$
y_1 y_2\,
T^{ 4+3R/2} 
\left(
{
y_{1}\over T
}
\right)^{\cC_{1}}
\left(
{y_{2}
\over
T}
\right)^{\gk'_{1}}
T^{\vep}
.
$$
Again, taking $\gk_{2}=0$ and $\gk_{1}>0$, the second inequalities in \eqref{eqs:gkToCs2}  are satisfied, and we can take $\gk_{1}=\cC_{2}/2<1$, so that $\gk_{1}'=\cC_{2}$. There are no new constraints on $R$.

Hence we see that the contribution by $\cR_{1}$ is dominated by that of $\fM$. The same holds for $\cR_{2}$ and $\cR_{3}$ by symmetry, and we have established the following crucial bound.
\begin{prop}
Fix $R\ge 10$ and $\vep > 0.$
For
any $y_{1},y_{2}>0$, $T\gg 1$,
and 
any
$0<\cC_{1},\cC_{2}<2$,
we have
\be\label{eq:range2}%ann
p_{_{T,R}}(y_{1},y_{2})
%&
\ll_{\cC_{1},\cC_{2},R,\vep}
%&
%\boxed{
%\boxed{
y_1y_2\,
 T^{9/2+3R/2}
\left(
y_{1}\over T
\right)^{\cC_{1}}
\left(
y_{2}\over T
\right)^{\cC_{2}}
T^{\vep}.
%}}
\ee%ann
In particular, $p_{T,R}$ satisfies  \eqref{eq:pBndAt0}, as needed in the Kuznetsov formula.
\end{prop}

\begin{rmk}
One can pull further and analyze contributions from higher poles. In so doing, the power in $T$  increases, so that the residual contributions $\cR$ dominate the contribution from $\fM$. It may still be possible to get further improvements from such an analysis, but the above is sufficient for our purposes, so we stop here.
\end{rmk}

%%%%%%%%%%%%%%%%%%%%%%%%%%%%%%%%%%%%%%%%
%%%%%%%%%%%%%%%%%%%%%%%%%%%%%%%%%%%%%%%%

\newpage

\section{Bounds for the Kloosterman and Eisenstein Contributions}\label{sec:boundsKloos}

Since we showed in the previous section that our choice $p_{T,R}$ of test function satisfies the requisite bound  \eqref{eq:pBndAt0}, we 
now invoke Kuznetsov's formula with this choice, and estimate the resulting components.

\subsection{Bounds for the Kloosterman integrals $\cJ$ and $\tilde \cJ$}\

We shall apply  the estimates obtained in the previous section to bound the Kloosterman integrals  $\widetilde\cJ$ and $\cJ$ defined in \S\ref{sec:preKuz}. 
%In our application $p=q=p_{_{T,R}}$ is determined by \ref{eq:pSharpIs}.  
%We write $\cJ_{p,p} =\cJ_{T,R}$ and $\widetilde\cJ_{p,p} = \widetilde\cJ_{T,R}$. There will be two bounds obtained which give better results depending on whether the arguments are small or large, but both bounds hold in all cases.
%
%\vskip 10pt\noindent
%$\underline{\text{{\bf Good bounds for   $\cJ_{T,R}(A_1, A_2)$ and $\widetilde\cJ_{T,R}(A)$ when   $A_1A_2$ (resp. $A$) is small.}}}$ 
%
%
%\vskip 10pt
%The more difficult case is $\cJ_{p,p}$ which we begin with.
We begin with an analysis of the more difficult case of $\cJ$. 
For $\gep_{1},\gep_{2}\in\pm1$, recall $\cJ_{\gep_{1},\gep_{2}}$ is given by
\bea
\nonumber
\cJ_{\gep_{1},\gep_{2}}(A_{1},A_{2})%; \gep_{1},\gep_{2})
&:=&
(A_{1}A_{2})^{-2}
\iint_{\R_{+}^{2}}
\iiint_{\R^{3}}
e(-\gep_{1}A_{1}x_{1}y_{1}-\gep_{2}A_{2}x_{2}y_{2})\,
\overline{p_{_{T,R}}}(A_{1}y_{1},A_{2}y_{2})
\\
\nonumber
&&\times
e\left(
-{A_{2}\over y_{2}}
{
x_{1}x_{3}+x_{2}
\over
x_{3}^{2}+x_{2}^{2}+1
}
%\right)
%e\left(
-{A_{1}\over y_{1}}
{
x_{2}(x_{1}x_{2}-x_{3})+x_{1}
\over
(x_{1}x_{2}-x_{3})^{2}+x_{1}^{2}+1
}
\right)
\\
\nonumber
&&\times
p_{_{T,R}}\left(
{A_{2}\over y_{2}}
{
\sqrt{(x_{1}x_{2}-x_{3})^{2}+x_{1}^{2}+1}
\over
x_{3}^{2}+x_{2}^{2}+1
}
,
{A_{1}\over y_{1}}
{
\sqrt{x_{3}^{2}+x_{2}^{2}+1}
\over
(x_{1}x_{2}-x_{3})^{2}+x_{1}^{2}+1
}
\right)
dx{dy\over y}
.
\\
\label{eq:Jis}
\eea
Here $dx = dx_1 dx_2 dx_3.$

We put absolute values inside the integral and note that the resulting bounds are then independent of $\gep_{1},\gep_{2}$, so it is convenient to drop $\gep_{1},\gep_{2}$ from the notation.
Since $p=p_{T,R}$, it is convenient to recall the dependence of $\cJ$ on $T$ and $R$, so we write henceforth 
$\cJ_{T,R}(A_1,A_2),$ etc. 
We have
$$
|\cJ_{T,R}(A_1, A_2)|\le 
\frac{1}{(A_{1}A_{2})^{2}}\iint_{\R_{+}^{2}}
\iiint_{\R^{3}} |p_{_{T,R}}(A_{1}y_{1},A_{2}y_{2})|\left|p_{_{T,R}}\left({A_{2}\over y_{2}}{\xi_{1}^{1/2}\over \xi_{2}},{A_{1}\over y_{1}}{\xi_{2}^{1/2}\over \xi_{1}}\right)\right|
dx{dy_{1}dy_{2}\over y_{1}y_{2}}
,
$$
where 
$$
\xi_{1}=
1
+x_{1}^{2}
+(x_{1}x_{2}-x_{3})^{2}
,
\qquad
\xi_{2}=
1
+x_{2}^{2}
+x_{3}^{2}
.
$$

For $i = 1,2$, break the $y_i$ integrals according to $y_i>1$ or $y_i<1$; this gives
$$
\cJ_{T,R}\le \cJ_{0}+\cJ_{1}+\cJ_{1}'+\cJ_{2},
$$
where the $y$ integral in $\cJ_{0}$ has $y_{1},y_{2}<1$, the $\cJ_{1},\cJ_{1}'$ integrals have one bigger and one smaller, and $\cJ_{2}$ has $y_{1},y_{2}>1$.
\\

%\subsubsection{Bounding $\cJ_{0}$}
%5\

We first estimate $\cJ_{0}$. 
Apply the bound in \eqref{eq:range2}
 to the second appearance of
 $p_{_{T,R}}$,
 choosing the largest possible values $\cC_{1}=\cC_{2}=
%\eta\in(0,2)
2-\vep
$:
\beann
\cJ_{0}
&\ll_{\vep}& 
{T^{9/2+3R/2+\vep}\over(A_{1}A_{2})^{2}}
\left(
{A_{1}
A_{2}\over T^{2}}
\right)^{2-\vep}
\int_{0}^{1}\int_{0}^{1} 
|p_{_{T,R}}(A_{1}y_{1},A_{2}y_{2})|
{A_{1}A_{2}\over (y_{1}y_{2})^{3-\vep}}
\iiint 
%
%\left|
%
{1\over (\xi_{1}\xi_{2})^{\frac{3-\vep}{2}}}
%
%
%p\left({A_{2}\over y_{2}}{\xi_{1}^{1/2}\over \xi_{2}},{A_{1}\over y_{1}}{\xi_{2}^{1/2}\over \xi_{1}}\right)
%
%\right|
%
dx
{dy_{1}dy_{2}\over y_{1}y_{2}}
\\
&\ll& 
T^{1/2+3R/2+3\vep}
\left(
A_{1}
A_{2}
\right)^{1-\vep}
\int_{0}^{1}\int_{0}^{1} 
{|p_{_{T,R}}(A_{1}y_{1},A_{2}y_{2})|
\over (y_{1}y_{2})^{3-\vep}}
{dy_{1}dy_{2}\over y_{1}y_{2}}
,
\eeann
since the $x$ integral converges absolutely.
Now estimate the first $p_{_{T,R}}$ again using \eqref{eq:range2}, with $\cC_{1}=\cC_{2}=2-\vep/2$:
\beann
\cJ_{0}
&\ll& 
{T^{1/2+3R+3\vep}}
\left(
{A_{1}}
{A_{2}}
\right)^{1-\vep}
\int_{0}^{1}\int_{0}^{1} 
A_{1}A_{2}y_{1}y_{2}
 T^{9/2+3R/2+\vep}
\left(
{A_{1}A_{2}y_{1}y_{2}\over T^{2}}
\right)^{2-\vep/2}
%
%
%|p(A_{1}y_{1},A_{2}y_{2})|
%
{
1\over (y_{1}y_{2})^{3-\vep}}
{dy_{1}dy_{2}\over y_{1}y_{2}}
\\
&\ll& 
{T^{1+3R+\vep}}
\left(
{A_{1}}
{A_{2}}
\right)^{4-\vep}
\int_{0}^{1}\int_{0}^{1} 
 (y_{1}y_{2})^{\vep/2}
{dy_{1}dy_{2}\over y_{1}y_{2}}
\ll
{T^{1+3R+\vep}}
\left(
{A_{1}}
{A_{2}}
\right)^{4-\vep}
,
\eeann
since the $y$-integral converges absolutely.

To bound $\cJ_{1}, \cJ_{1}',$ and $\cJ_{2}$, 
we simply follow the same procedure as above with minimal changes to ensure convergence, as follows.
For $\cJ_{2}$, in the second application of \eqref{eq:range2}, we choose $\cC_{1}=\cC_{2}=2-2\vep$, say, so the final $y$-integral converges absolutely.
Similarly, for $\cJ_{1}$ and $\cJ_{1}'$,  we do the same as before, except in the second application of \eqref{eq:range2}, we choose $\cC_{1}=2-\vep/2$, $\cC_{2}=2-2\vep$ (or vice versa), so that the $y$-integral converges absolutely. 
%
%
%\begin{comment}
%Taking $\eta,\eta'$ as large as possible (which is $2$), we obtain
We have thus proved that
%\begin{prop}
%For  $T\gg1$, $R\ge10$, any $A_{1},A_{2}>0$, and any $\vep>0$ sufficiently small, we have
$$
\cJ_{T,R}(A_{1},A_{2}) \;
\ll_{R,\vep} \;
T^{1+3R+\vep}
(A_{1}A_{2})^{4-\vep}.$$
%\end{prop}
%\end{comment}

\vskip 10pt
Next, we want a similar bound  for $\widetilde\cJ_{T,R}(A).$
Recall that we have
\begin{align*}
\widetilde\cJ_{T,R}
(A) & \;
\le \;
A^{-2}
\iint\limits_{\R_{+}^{2}}
\iint\limits_{\R^{2}}
|p_{_{T,R}}(Ay_{1},y_{2})| \\
&
\hskip 30pt \times
\left|p_{_{T,R}}\left(
{y_{2}}\cdot{
{\sqrt{1 +  x_{1}^{2}+x_{2}^{2} }}\over 1+ x_{1}^{2}}
, \;\;
\frac{A}{y_{1}y_{2}}\cdot{
\sqrt{1 + x_{1}^{2}}\over 1+ x_{1}^{2}+x_{2}^{2}}
\right)\right| 
dx_{1}dx_{2}{dy_{1}dy_{2}\over y_{1}y_{2}^{2}}.
\end{align*}
Note that here the integral involves $dy_{2}/y_{2}^{2}$, whereas in $\cJ$ the integral has $dy_{2}/y_{2}$. This will result in a weaker final bound.

As before, for $i = 1,2$, break the $y_i$ integrals according to $y_i>1$ or $y_i<1$; this gives
$$
\widetilde \cJ_{T,R}\le \widetilde \cJ_{0}+\widetilde\cJ_{1}+\widetilde\cJ_{1}'+\widetilde\cJ_{2},
$$
where the $y$ integral in $\widetilde\cJ_{0}$ has $y_{1},y_{2}<1$, the $\widetilde\cJ_{1},\widetilde\cJ_{1}'$ integrals have one bigger and one smaller, and $\widetilde\cJ_{2}$ has $y_{1},y_{2}>1$.
\\

%\subsubsection{Bounding $\cJ_{0}$}
%5\

We first bound $\widetilde\cJ_{0}$. Set $\xi_1 = 1+x_1^2$ and $\xi_2 = 1+x_1^2+x_2^2.$
Replace the second $p_{_{T,R}}$ by its bound in \eqref{eq:range2}, with the choice $\cC_{1}=\vep$ and $ \cC_{2}=2-\vep
%\eta\in(1,2)
$:
\beann
\widetilde\cJ_{0}
&\ll_{\vep}& 
\hskip -5pt
{T^{9/2+3R/2+\vep}\over A^{2}}
\int_{0}^{1}\int_{0}^{1} 
\int\limits_{-\infty}^\infty \int\limits_{-\infty}^\infty 
|p_{_{T,R}}(Ay_{1}, y_{2})| 
\frac{A}{y_1 \sqrt{\xi_1\xi_2}} 
\left(\frac{y_2\sqrt{\xi_2}}{T\xi_1}  \right)^\vep
\left(
\frac{A\sqrt{\xi_1}}{y_1y_2\xi_2T}   
\right)^{2-\vep}  dx_1 dx_2
{dy_{1}dy_{2}\over y_{1}y_{2}^2}
\\
&\ll& \hskip -5pt
{T^{5/2+3R/2+\vep}}
{A^{1-\vep}}
\int_{0}^{1}\int_{0}^{1} 
{|p_{_{T,R}}(Ay_{1}, y_{2})|\,\over
 (y_{1}y_2)^{3-\vep}}\;
{dy_{1}dy_{2}\over y_{1}y_{2}}
,
\eeann
since again the $x$ integral converges absolutely. Here we used that for fixed $x_{1}>1$ and $Z<-1/2$,
$$
\int_{\R}(1+x_{1}^{2}+x_{2}^{2})^{Z}dx_{2}\ll x_{1}^{1+2Z}.
$$

Now apply \eqref{eq:range2} to $p_{_{T,R}}(Ay_1, y_2)$ in the above integral with $\cC_1 = \cC_2 = 2-\vep/2$.
 %with $\eta' \in  (\eta, 2).$ 
 It follows that
$$
\widetilde\cJ_{0} \; 
\ll 
\; 
T^{5/2+3R+\vep} 
A^{1-\vep} 
\int_0^1\int_0^1 
Ay_{1}y_{2}T^{9/2+3R/2+\vep}
\left(
{
Ay_{1}y_{2}
\over
T^{2}
}
\right)^{2-\vep/2}
(y_{1}y_{2})^{-3+\vep}
\frac{dy_1 dy_2}{y_1y_2}.
$$
%Since $\eta' > \eta$, t
Since $\vep>0$, the above $y$-integral converges, and we obtain the bound 
$$
\widetilde\cJ_{0} \; \ll \; T^{3+3R+\vep} A^{4-\vep}.
$$

Then we  bound $\widetilde\cJ_{1}, \widetilde\cJ_{1}',$ and $\widetilde\cJ_{2}$, 
%we simply follow 
by
the same procedure as above, with suitable modifications, as before.
%For $\cJ_{2}$,  we just choose $\eta'<\eta$, so the final integral converges absolutely.
%
%For $\widetilde\cJ_{1}$ and $\widetilde\cJ_{1}'$,  we do the same as before, except in the second application of \eqref{eq:range2}, we choose $\eta'_{1}>\eta$, $\eta_{2}'<\eta$ (or vice versa), so that the integral converges absolutely. 
%\vskip 10pt
  We record the previous computations in the following.

\begin{proposition}\label{Bnd1KloostIntegral} Fix $R\ge 10$, and any small $\vep > 0$.
\vskip 4pt
For any $A_{1},A_{2}>0$ and $T\gg 1$,, we have
\be
\cJ_{T,R}(A_{1},A_{2})
\ll_{R,\vep}
T^{1+3R+\vep}
(A_{1}A_{2})^{4-\vep}
.
\ee

For any $A> 0$ and  $T\gg 1$, we have
\be\label{eq:tildecJbnd1}
\widetilde\cJ_{T,R}(A)
\ll_{R,\vep}
T^{3+3R+\vep}
A^{4-\vep}.
\ee
\end{proposition}

\subsection{Bounds for the Kloosterman contributions}\label{subs:pfKuz}\

\vskip 10pt

We shall now apply the bounds
  in Proposition  \ref{Bnd1KloostIntegral} %,   \ref{Bnd2KloostIntegral},
 to estimate the Kloosterman contributions $\cK, \tilde \cK,$ and $\tilde \cK^{\scriptscriptstyle \vee}$
\begin{comment}
\bea
\label{Kloos}
\widetilde\cK
&=&
\sum_{\gep=\pm1}
\sum_{D_{1}\mid D_{2}\atop m_{2}D_{1}^{2}=n_{1}D_{2}}
{
\widetilde S(\gep m_{1},n_{1},n_{2},D_{1},D_{2})
\over
D_{1}D_{2}
}
\widetilde\cJ_{p,q}\left(\sqrt{n_{1}n_{2}m_{1}\over D_{1}D_{2}}; \, \gep\right),
\\
\nonumber
\widetilde\cK^{\scriptscriptstyle \vee}
&=&
\sum_{\gep=\pm1}
\sum_{D_{2}\mid D_{1}\atop m_{1}D_{2}^{2}=n_{2}D_{1}}
{
\widetilde S(\gep m_{2},n_{2},n_{1},D_{2},D_{1})
\over
D_{1}D_{2}
}
\widetilde\cJ_{p, q}^{\scriptscriptstyle \vee}\left(\sqrt{n_{1}n_{2}m_{2}\over D_{1}D_{2}};\, \gep\right),
\\
\nonumber
\cK
&=&
\sum_{\gep_{1},\gep_{2}=\pm1}
\sum_{D_{1},D_{2}}
{
S(\gep_{1}m_{1},\gep_{2}m_{2},n_{1},n_{2},D_{1},D_{2})
\over
D_{1}D_{2}
}
\cJ_{p, q}
\left(
{
\sqrt{m_{1}n_{2}D_{1}}
\over
D_{2}
}
,
{
\sqrt{m_{2}n_{1}D_{2}}
\over
D_{1}
};\, \gep_{1},\gep_{2}
\right)
,
\eea
\end{comment}
appearing in the geometric side of the $\GL(3)$ Kuznetsov formula \eqref{eq:cMis}. We need estimates  for these contributions using our choice of test function (or its transform)
%$p = q = p_{_{T,R}}$ with $p_{T,R}^\#$ 
given by \eqref{eq:pSharpIs}.

 Let us begin by bounding the long element Kloosterman contribution $\cK.$ %The simplest method uses
We only use  the trivial bound for the Kloosterman sum:
\begin{equation} \label{TrivialBound}
|S(m_1,n_1, m_2,n_2,D_1, D_2)|
\ll_\vep
(D_1D_2)^{1+\vep}
.
\end{equation}

  It immediately follows from Proposition \ref{Bnd1KloostIntegral} and the trivial bound \eqref{TrivialBound}
 that 
  
 \begin{align}\cK & \ll \sum_{D_{1}=1}^\infty
\sum_{D_{2}=1}^\infty
\frac{|S(m_1,m_2, n_1,n_2,D_1, D_2)|}{D_1D_2} \left|\cJ_{T,R}
\left(
{\sqrt{m_1n_2D_{1}}
\over D_{2}}
,
{\sqrt{m_2n_1D_{2}}
\over D_{1}}
\right)
\right| \notag \\
& \ll_{R,\vep}  |m_1m_2n_1n_2|^2 \, T^{1+3R+\vep} \sum_{D_{1}=1}^\infty
\sum_{D_{2}=1}^\infty |D_1D_2|^{\vep -2}\notag \\\label{eq:KloostBnd1}
& \ll |m_1m_2n_1n_2|^2 \,T^{1+3R+\vep}. 
\end{align}

Next, we obtain a similar proposition for the lower rank Kloosterman contributions $\widetilde{\cK}$ and $\widetilde\cK^{\scriptscriptstyle \vee}$. In this case, we have the trivial bound
$$\widetilde S(m_1,n_1,n_2, D_1, D_2) \ll_\vep (D_1D_2)^{1+\vep}.$$
 It immediately follows from Proposition \ref{Bnd1KloostIntegral} that 
  
 \begin{align}
 \widetilde\cK & \ll \sum_{D_2=1}^\infty\sum_{D_{1}\mid D_{2}\atop m_{2}D_{1}^{2}=n_{1}D_{2}}
{
|\widetilde S(m_{1},n_{1},n_{2},D_{1},D_{2})|
\over
D_{1}D_{2}
}
\left|
\widetilde\cJ_{T,R}\left(\sqrt{n_{1}n_{2}m_{1}\over D_{1}D_{2}}\right)
\right|
 \notag \\
& \ll_{R,\vep}  |m_1n_1n_2|^2 \, T^{3+3R+\epsilon} \sum_{D_2=1}^\infty\sum_{D_{1}\mid D_{2}\atop m_{2}D_{1}^{2}=n_{1}D_{2}}
 |D_1D_2|^{2\epsilon -2}\notag \\\label{eq:KloostBnd1}
& \ll |m_1n_1n_2|^2 \,T^{3+3R+\epsilon}. 
\end{align}

\begin{comment}
To get a better bound in the $|m_1n_1n_2|$ aspect, we break the sum for $\widetilde\cK$ into two pieces
$\widetilde\cK = \widetilde \cK_0 + \widetilde\cK_1$ with
$$\widetilde \cK_0 = \sum_{D_2 < H}, \qquad \widetilde\cK_1 = \sum_{D_2 \ge H}.$$
On the first piece we apply Proposition  \ref{Bnd2KloostIntegral} to obtain
$$\widetilde\cK_0 \; \ll \; T^{9+3R+\vep} \sum_{D_2 < H} \sum_{D_1 |D_2} (D_1D_2)^\vep \; \ll \; T^{9+3R+\vep} \, H^{1+\vep}.$$
On the second piece, we apply Proposition \ref{Bnd1KloostIntegral} to obtain 
$$\widetilde\cK_1 \; \ll \; |m_1n_1n_2|^2 \, T^{3+3R+\vep} \sum_{D_2 \ge H} \sum_{D_1 |D_2} |D_1D_2|^{2\epsilon -2} \; \ll \; |m_1n_1n_2|^2 \, T^{3+3R+\vep} H^{-1+\vep}.$$
Choosing $H = |m_1n_1n_2| T^{-3}$ gives the bound
$$\widetilde\cK \; \ll \; |m_1n_1n_2|^{1+\vep} T^{6+3R+\vep}.$$

\end{comment}

The same argument applies to $\widetilde\cK^{\scriptscriptstyle \vee}$.
 We have proved the following.
 \begin{proposition}\label{Bnd2KloostContr} Fix $R\ge 10, \, T \gg 1,$ and $\vep > 0.$ 
Then the
 %lower rank
  Kloosterman contributions $\cK$, $\widetilde\cK
   %= \widetilde\cK_{T,R}(m_1,n_1,n_2)
   $ and 
   $\widetilde\cK^{\scriptscriptstyle \vee}%_{T,R}(m_2,n_2,n_1)
   $
   %, as defined 
   in
 \eqref{eq:cMis} 
  satisfy the bounds:
\beann
\cK
&\ll_{R,\vep}&
|m_1m_2n_1n_2|^2 \,T^{1+3R+\vep}
,
\\
 \widetilde\cK %_{T,R}(m_1,n_1,n_2) 
& \ll_{R,\vep}& |m_1n_1n_2|^2\, T^{3+3R+\vep}, 
\\
 \cK^{\scriptscriptstyle \vee}
 %_{T,R}(m_2,n_2,n_1) 
 &\ll_{R,\vep} 
 &
 |m_2n_2n_1|^2 \,T^{3+3R+\vep}
 .
 \eeann
% $$\widetilde\cK_{T,R}(m_2,n_2,n_1) \ll_{\vep}  |m_1n_1n_2|^{1+\vep} \,T^{6+3R+\vep}, \qquad \cK^{\scriptscriptstyle \vee}_{T,R}(m_2,n_2,n_1) \ll_{\vep}  |m_2n_2n_1|^{1+\vep}\, T^{6+3R+\vep}
%.$$
\end{proposition}

%%%%%%%%%%%%%%%%%%%%%%%%%%%%%%%%%%%%%%%%
%%%%%%%%%%%%%%%%%%%%%%%%%%%%%%%%%%%%%%%%

\subsection{Bounds on the Continuous Spectrum}\label{subs:BndsCont}\

\vskip 10pt
 %In this section 
 Next we obtain bounds for the terms $\cE_{min},\, \cE_{max}$ coming from the continuous spectrum in the Kuznetsov  formula \eqref{eq:Kuz}.
 % with the choice of test function
%$p = q = p_{_{T,R}}$ with $p_{T,R}^\#$ given by \eqref{eq:pSharpIs}. 
%Let's 
We begin with the term coming from the minimal parabolic Eisenstein series:
\be\label{eq:cEmin}
\cE_{min}
=
{1%c_{1}
\over (4\pi i)^{2}}
\int\limits_{-i\infty}^{i\infty}
\int\limits_{-i\infty}^{i\infty}
A_{\nu}(n_{1},n_{2})
\overline{A_{\nu}(m_{1},m_{2})}
{
|p_{T,R}^{\sharp}(\nu_{1},\nu_{2})|^{2}
\over \prod\limits_{k=1}^{3}
\left|
\gz(1+3\nu_{k})
\G\left({1+3\nu_{k}\over 2}\right)
\right|^{2}
%\prod_{k=1}^{3}
%
%\G\left({-3\nu_{k}^{(j)}\over 2}\right)
}\;
d\nu_{1}
d\nu_{2}
,
\ee
where % $c_{1}$ is an absolute constant, and 
the minimal Eisenstein coefficients satisfy 
\eqref{eq:AnuBnd}.
\begin{comment}
$$
|A_{\nu}(n_{1},n_{2})|
\ll_{\vep}
(n_{1}n_{2})^{\vep}
.
$$
\end{comment}
%
%Recall from
Inserting  the choice of test function \eqref{eq:pSharpIs}
into \eqref{eq:cEmin}
%that
%\be
%p_{T,R}^{\sharp}(\nu_{1},\nu_{2})
%:=
%\sqrt6\
%e^{\ga_{1}^{2}+\ga_{2}^{2}+\ga_3^2\over 2T^{2}}
%\prod_{1\le j\le 3} 
%\G\left({2+R+3\nu_j\over4}\right)
%\G\left({2+R-3\nu_j\over4}\right)
%.\notag
%\ee
%
%
%It follows from Stirling's  formula , together with
and using
 the de la Vall\'ee Poussin bound (Prime Number Theorem)
$$|\zeta(1+it)| \gg \frac{1}{ \log(2+|t|)},$$
we get from Stirling's formula \eqref{Stirling} that
\begin{align}
\cE_{min}
& \ll
{1%c_{1}
\over (4\pi i)^{2}}
\int\limits_{-i T^{1+\vep}}^{i T^{1+\vep}}
\int\limits_{-i T^{1+\vep}}^{i T^{1+\vep}}
|m_1m_2n_1n_2|^\vep \;\frac{\prod\limits_{j=1}^3 \left|\G\left({2+R+3\nu_j\over4}\right)
\G\left({2+R-3\nu_j\over4}\right)\right|^2
}{\prod\limits_{k=1}^{3}
\left|
\gz(1+3\nu_{k})
\G\left({1+3\nu_{k}\over 2}\right)
\right|^{2}} \; d\nu_1 d\nu_2\notag\\
& \ll T^{2+3R+\vep} \, |m_1m_2n_1n_2|^\vep. \label{Emin}
\end{align}

%\vskip 5pt
Next we consider the term coming from the maximal parabolic Eisenstein series:
\beann
\cE_{max}
&=&
{c\over 2\pi i}
\sum_{j=1}^\infty\;
\int\limits_{-i\infty}^{i\infty}%\nu\in  i\R}
{B_{\nu,r_{j}}(n_{1},n_{2})
\overline{B_{\nu,r_{j}}(m_{1},m_{2})}
\over
L(1, \Ad u_{j})
|L(1+3\nu, u_{j})|^{2}
}\;
{
\left|p^{\sharp}_{T,R}\left(\nu- \frac{ir_{j}}{3}, \; \frac{2ir_{j}}{3}\right)\right|^{2}
\over 
\left|
\G\left({1+3\nu-ir_{j}\over2}\right)
\G\left({1+2ir_{j}\over2}\right)
\G\left({1+3\nu+ir_{j}\over2}\right)
\right|^{2}
}\;
d\nu
,
\eeann
where $c$ is an  absolute constant,  $\{u_{j}\}$ is a basis of Hecke-Maass forms for $\GL(2,\Z)$ (each of eigenvalue $1/4+r_{j}^{2}$), and
the Fourier coefficients satisfy
\eqref{eq:BnuBnd}.
%$$
%|B_{\nu,r_{j}}(n_{1},n_{2})|
%\ll
%(n_{1}n_{2})^{\vep}.
%$$
Here we have the lower bounds
$$L(1, \Ad u_{j}) \gg_\vep (1+|r_j|)^{-\vep},\qquad
L(1+3\nu, u_j) \gg_\vep (1 + |\nu| + |r_j|)^{-\vep}.$$
These lower bounds follow from
 \cite{HoffsteinLockhart1994, HoffsteinRamakrishnan1995, JacquetShalika1976} and \cite{GelbartLapidSarnak2004}.
Combining the above lower bounds with Stirling's formula \eqref{Stirling}, it follows that
\begin{align}
\cE_{max}
& \ll  |m_1m_2n_1n_2|^{1/2+\vep} \, \sum_{r_j \ll T^{1+\vep}} \int\limits_{-iT^{1+\vep}}^{iT^{1+\vep}} \frac{ \left|\Gamma\left( \frac{2+R +3\nu-ir_j}{4} \right)\right|^8 
}{L(1, \Ad u_{j})
|L(1+3\nu, u_{j})|^{2}}\notag\\
& \hskip 110pt\times \frac{\left|\Gamma\left( \frac{2+R + 2ir_j}{4} \right)\right|^4 }{\left|
\G\left({1+3\nu-ir_{j}\over2}\right)
\G\left({1+2ir_{j}\over2}\right)
\G\left({1+3\nu+ir_{j}\over2}\right)
\right|^{2}} \; d\nu\notag\\ 
& \notag\\
& \ll |m_1m_2n_1n_2|^{1/2+\vep} \, \sum_{r_j \ll T^{1+\vep}} \int\limits_{-iT^{1+\vep}}^{iT^{1+\vep}} \, \big(1+\left|3\nu - ir_j   \right|\big)^{\frac{R}{4}\cdot8} \cdot \big(1 +|r_j|  \big)^{\frac{R}{4}\cdot 4} \; |d\nu|
\notag\\
&
\notag\\
& \ll |m_1m_2n_1n_2|^{1/2+\vep} \; T^{3+3R+\vep}\label{Emax}
,\end{align}
using  Weyl's Law for $\GL(2)$ and the Ramanujan conjectures at infinity (for $\GL(2)$). In summary, we have proved the following.
%\vskip 10pt
\begin{proposition} \label{EisContr} Fix $R\ge 10$ and $\vep > 0.$ For any $T\gg 1$, we have
 $$\cE_{min} \ll_{R,\vep}  |m_1m_2n_1n_2|^\vep\, T^{2+3R+\vep}, \qquad \cE_{max} \ll_{R,\vep}\,  |m_1m_2n_1n_2|^{1/2+\vep} T^{3+3R+\vep}.$$
\end{proposition}

%%%%%%%%%%%%%%%%%%%%%%%%%%%%%%%%%%%%%%%%
%%%%%%%%%%%%%%%%%%%%%%%%%%%%%%%%%%%%%%%%

%%%%%%%%%%%%%%%%%%%%%%%%%%%%%%%%%%%%%%%%
%%%%%%%%%%%%%%%%%%%%%%%%%%%%%%%%%%%%%%%%

\newpage

\section{Proof of Theorem \ref{thm:Kuz}}\label{sec:PfThm1}

\subsection{The Main Term}\
\vskip 10pt
  We begin by computing an asymptotic  formula for the main term $\cM$ in the Kuznetsov  formula \eqref{eq:Kuz}, \eqref{eq:cMis}. 
  It follows from \eqref{eq:KLParseval} and Stirling's asymptotic formula  \eqref{Stirling} that  the inner product in the main term \eqref{eq:cMis}
 becomes
\begin{align} \label{eq:MainTerm}
\<p_{_{T,R}}, \, p_{_{T,R}}\>
&= \frac{1}{(\pi i)^2}
\int\limits_{-i\infty}^{i\infty} \int\limits_{-i\infty}^{i\infty}
|p_{T,R}^{\sharp}(\nu_{1},\nu_{2})|^{2}
{d\nu_{1}d\nu_{2}
\over
\prod\limits_{1\le j\le 3}
\G\left(
{3\nu_j\over 2}
\right)
\G\left(
{-3\nu_j\over 2}
\right)
}  \; d\nu_1\nu_2
\\
&\sim
c'
\iint\limits_{|\nu_{1}|, \,|\nu_{2}|\ll T}
{
\left|
\prod\limits_{1\le j\le 3} 
\G\left({2+R+3\nu_j\over4}\right)
\G\left({2+R-3\nu_j\over4}\right)
\right|^{2}
\over
\prod\limits_{1\le j\le 3}
\G\left(
{3\nu_j\over 2}
\right)
\G\left(
{-3\nu_j\over 2}
\right)
}
d\nu_{1}d\nu_{2}\notag\\
&\notag \\
& \sim
c\,
T^{5+3R}
,\notag
\end{align}
for some constants $c, c' > 0.$ This gives the  $T^{5+3R}$ main term as claimed in \eqref{eq:hTmain}.\

\subsection{Completion of the proof of Theorem 1.2}\
\vskip 10pt

It follows from the Kuznetsov  formula 
\eqref{eq:Kuz}, together with the choice of test function
$p = %q = 
p_{_{T,R}}$ (with $p_{T,R}^\#$ given by \eqref{eq:pSharpIs}), that
\begin{align} \label{Thm1.2proof}
\cC & =
\sum_{j}
{A_{j}(m_{1},m_{2})
\overline{A_{j}(n_{1},n_{2})}
}
{h_{T,R}(\nu_{1}^{(j)},\nu_{2}^{(j)})\over \cL_{j}}\notag\\
& = \cM
+
\cK
+
\widetilde{\cK}+
\widetilde\cK^{\scriptscriptstyle \vee} -\cE_{min}-\cE_{max}
,
\end{align}
with $\cM = \bo_{\left\{n_{1}=m_{1}\atop n_{2}=m_{2}\right\}}
 \<p_{_{T,R}}, \, p_{_{T,R}}\>$.

\vskip 5pt
Then % proof of
 Theorem  \ref{thm:Kuz} is an immediate consequence of the estimates
\begin{align*} 
 \<p_{_{T,R}}, \, p_{_{T,R}}\> &  \sim c T^{5+3R}\\
 |\cE_{min}|, |\cE_{max}| & \ll_{\vep}|m_1m_2n_1n_2|^{1/2+\vep}\ T^{3+3R+\vep}\\
 \big|\cK\big| + \big|\widetilde{\cK}\big| + \big|\widetilde\cK^{\scriptscriptstyle \vee}\big|  & \ll_{R,\vep} |m_1m_2n_1n_2|^2 \, T^{3+3R+\vep}\end{align*}
 
 \noindent
given in \eqref{eq:MainTerm}, and  Propositions  \ref{EisContr},  %\ref{Bnd1KloostContr}, 
\ref{Bnd2KloostContr}, respectively.

%%%%%%%%%%%%%%%%%%%%%%%%%%%%%%%%%%%%%%%%
%%%%%%%%%%%%%%%%%%%%%%%%%%%%%%%%%%%%%%%%

\newpage

\section{The Explicit Formula}\label{subs:ExpForm}

%%%%%%%%%%%%%%%%%%%%%%%%%%%%%%%%%%%%%%%%
%%%%%%%%%%%%%%%%%%%%%%%%%%%%%%%%%%%%%%%%

%\newpage

For a Hecke-Maass form $\phi$ on $\GL(3)$, let $\rho(\phi)$ be one of
$$
\rho(\phi)
=
\threecase
{\phi}{}
{\sym^{2}\phi}{}
{\Ad\phi,}{}
%{\phi\times\phi,}{}
$$
and let $L(s,\rho(\phi))$ be the corresponding $L$-function. Note that the  Maass form dual to $\phi$ is just the complex conjugate $\bar\phi.$ In view of the identity
$$L(s, \Ad\phi) = \frac{L(s, \,\phi\times \bar\phi)}{\zeta(s)},$$
it is easier to  work with the Rankin-Selberg convolution of $\phi$ and $\bar\phi$ instead of the adjoint L-function.

 Define
 $$\Lambda(s, \,\rho(\phi)) := \begin{cases} \pi^{-\frac{3s}{2}} \prod\limits_{k=1}^3 \Gamma\left(\frac{s+\alpha_k}{2}  \right)^{-1}  \prod\limits_p\prod\limits_{k=1}^3 \left(1 - \frac{\alpha_k(p)}{p^s}  \right)^{-1}, & {\rm if} \; \rho(\phi) = \phi,\\
 & \\
 \pi^{-\frac{5s}{2}} \prod\limits_{1\le j\le k\le 3} \Gamma\left(\frac{s+\alpha_j+\alpha_k}{2}  \right)^{-1}
  \prod\limits_p\prod\limits_{1\le j \le k \le 3} \left(1 - \frac{\alpha_j(p)\alpha_k(p)}{p^s}  \right)^{-1}, & {\rm if} \; \rho(\phi) = \sym^{2}\phi,
  \\
  & \\
  \pi^{-\frac{9s}{2}} \prod\limits_{j=1}^3 \, \prod\limits_{k=1}^3\G\left({s+\alpha_j- \alpha_k\over 2}\right)^{-1}
  \prod\limits_p \, \prod\limits_{j=1}^3 \, \prod\limits_{k=1}^3\left(1-{\ga_j(p)\overline{\ga}_k(p)\over p^s}\right)^{-1}, & {\rm if} \; \rho(\phi) = \phi\times \bar\phi.\end{cases}$$
  Then, in all the above cases, we have the functional equation
  $$\Lambda(s, \,\rho(\phi)) = \Lambda(1-s, \widetilde{\rho(\phi)}),$$
  where $\tilde \pi$ is the contragredient representation of $\pi$; its $L$-function has Dirichlet coefficients which are complex conjugates of the original.
This follows from \cite{GodementJacquet1972,  BumpGinzburg1992, JacquetPSShalika1983}, respectively.

We shall use the functional equation for $\Lambda(s, \,\rho(\phi))$  to determine the so-called ``explicit formula'' relating zeros and poles of $\Lambda(s, \,\rho(\phi))$ with sums over prime power Fourier coefficients of $L(s, \rho(\phi))$.
\vskip 5pt
Let $G$ be any holomorphic function in the region $-1\le \Re(s)\le2$ satisfying
$$
G(s)=G(1-s),\quad |s^2G(s)|\ll1.
$$
Let $\rho_i=\foh+i\g_i$ ($i =  \pm1, \pm2, \ldots$) run over the zeros of $\gL(s, \rho(\phi))$ with corresponding multiplicity. As we have assumed GRH, the ordinates form a real increasing sequence 
$$
\cdots\le\g_{-2}\le\g_{-1}\le0\le\g_1\le\g_2\le\cdots
$$
By the functional equation and standard shifts of contours, together with the fact (first proved by \cite{BumpGinzburg1992}) that $\frac{\gL'(s)}{\gL(s)}$ has at most simple poles at $s=0,1$, with residue 
\begin{equation} r_{\rho(\phi)} = \begin{cases} 0 & {\rm if} \;  \rho(\phi) = \phi,
\\
 0 & {\rm if} \;  \rho(\phi) = \sym^{2}\phi \; {\rm and} \; \phi \; {\rm not \; self \; dual,}\\
 1 & {\rm if} \; \rho(\phi) = \phi\times\bar\phi, \end{cases}\label{residues}
 \end{equation}
 we  have 
\begin{align}\label{eq1}
\sum_\rho G(\rho)- r_{\rho(\phi)}   \big(G(0) + G(1)\big) & = {1\over 2\pi i} \int\limits_{2-i\infty}^{2+i\infty}  G(s) \left[ {\gL'\over \gL}\big(s, \rho(\phi)\big) + {\gL'\over \gL}\big(s, \widetilde{\rho(\phi)}\big)\right] ds\\
&  = \sum_{p\le\infty}H_{\rho(\phi)}(p).
\end{align}
For a finite prime  $p < \infty,$ the function $H_{\rho(\phi)}(p)$ is defined by

$$H_{\rho(\phi)}(p) = \begin{cases}  -\sum\limits_{\ell=1}^\infty \left(\sum\limits_{k=1}^3 \left(\alpha_k(p)^\ell + \overline{\alpha_k(p)^\ell}\right) \right) F\left(p^\ell\right)) \log p,  & {\rm if}  \; \rho(\phi) = \phi,\\
& \\
-\sum\limits_{\ell=1}^\infty \left(\sum\limits_{1\le j \le k \le 3} \left(\alpha_j(p)^\ell \alpha_k(p)^\ell + \overline{\alpha_j(p)^\ell \alpha_k(p)^\ell}\right)\hskip -2pt\right)\hskip -2pt  F\left(p^\ell\right)) \log p, & {\rm if}  \; \rho(\phi) = \sym^{2}\phi,\\
& \\
-2\sum\limits_{\ell=1}^\infty \left( \sum\limits_{j=1}^3\sum\limits_{k=1}^3 \ga_j(p)^\ell \, \overline{\ga_k(p)^\ell} \right) F(p^\ell)\log p, & {\rm if}  \; \rho(\phi) = \phi\times \bar\phi.
\end{cases}$$

Here $F(y)$ is the inverse Mellin transform of $G(s)$,
$$
F(y)={1\over 2\pi i} \int\limits_{(1/2)} G(s) y^{-s} \;ds.
$$
For $p=\infty$, we have that $H_{\rho(\phi)}(\infty)$ equals
$$ = \begin{cases} -3 F(1) \log\pi + {1\over 4\pi i} \sum\limits_{k=1}^3  \;\int\limits_{\frac12-i\infty}^{\frac12+i\infty} \bigg({\G'\over \G}\left({s+\alpha_k\over 2}\right) +  {\G'\over \G}\left({s-\alpha_k\over 2}\right)\hskip -3pt\bigg) G(s)\, ds, & {\rm if}\; \rho(\phi) = \phi, 
\\
\-5 F(1) \log\pi + {1\over 4\pi i} \sum\limits_{1\le j\le k\le 3} \; \int\limits_{\frac12-i\infty}^{\frac12+i\infty} \bigg({\G'\over \G}\left({s+\alpha_j+\alpha_k\over 2}\right) +  {\G'\over \G}\left({s-\alpha_j-\alpha_k\over 2}\right)\hskip -3pt\bigg)G(s)\, ds, & {\rm if}\; \rho(\phi) = \sym^{2}\phi, \\
& \\
-9 F(1) \log\pi + {1\over 2\pi i} \sum\limits_{j=1}^3\sum\limits_{k=1}^3 \; \int\limits_{\frac12-i\infty}^{\frac12+i\infty} {\G'\over \G}\left({s+\alpha_j-\alpha_k\over 2}\right)G(s)\, ds, & {\rm if}\; \rho(\phi) = \phi\times \bar\phi.
\end{cases}
$$

Fix an even test function $\psi$ of Schwartz class whose Fourier transform has compact support, and apply the above formulae making the choice
$$
G(s) = \psi\left( \big(s-1/2\big) \,{\log C_{\rho}\over 2\pi i}\right)
,
$$
%where $R>1$ is a parameter to be specified later (roughly $R\approx c_f$)
where we recall (see \eqref{Arho}) that $C_{\rho}$ is the weighted average value of the conductor of $L(s, \rho(\phi)).$

Then
$$
F(y) = {1\over \sqrt y \log C_{\rho}%R
}\hat \psi\left({\log y\over \log C_{\rho}%R
}\right).
$$

It follows that 
$$H_{\rho(\phi)}(\infty) =\begin{cases} -\frac{3\widehat\psi(0) \log \pi}{\log C_{\rho}}  +  \sum\limits_{k=1}^3  \int\limits_{-\infty}^\infty \hskip -3pt\left(\frac{  {\G'\over \G}\left( {\pi i x\over \log C_{\rho}%R
} + {1+2\alpha_k\over 4} \right) +  {\G'\over \G}\left( {\pi i x\over \log C_{\rho}%R
} + {1-2\alpha_k\over 4} \right)}{ 2\log C_\rho }\right)\hskip -1pt
\psi(x) \, dx, & \rho(\phi) = \phi,\\
& \\
 -\frac{5\widehat\psi(0) \log \pi}{\log C_{\rho}} + \frac{1}{2\log C_\rho} \sum\limits_{1\le j\le k\le 3}\;  \int\limits_{-\infty}^\infty \bigg(    {\G'\over \G}\left( {\pi i x\over \log C_{\rho}%R
} + {1+2(\alpha_j + \alpha_k)\over 4} \right)\\
& 
\\
 \hskip 130pt+ {\G'\over \G}\left( {\pi i x\over \log C_{\rho}%R
} + {1-2(\alpha_j + \alpha_k)\over 4} \bigg)   \right) \hskip -1pt
\psi(x) \, dx, & \rho(\phi) = \sym^{2}\phi,\\
 & \\ -\frac{9\widehat\psi(0) \log \pi}{\log C_{\rho}} \; + \; \frac{1}{\log C_{\rho}}\sum\limits_{j=1}^3 \sum\limits_{k=1}^3  \int\limits_{-\infty}^\infty {\G'\over \G}\left( {\pi i x\over \log C_{\rho}%R
} + {1+2(\alpha_j - \alpha_k)\over 4} \right)
\psi(x) \, dx, & \rho(\phi) = \phi\times \bar\phi. \end{cases}$$

Recall that
\begin{align}
D(\rho(\phi); \,\psi):= \sum_{\Lambda\left(\frac12+i\gamma_i, \,\rho(\phi)\right) = 0}\psi\left(\gamma_i \,{\log C_{\rho}\over 2\pi} \right).
\end{align}
Consequently,  \eqref{eq1} becomes
\begin{align}\label{Dphipsi}
D(\rho(\phi); \,\psi) = \begin{cases} B_{\rho(\phi)} - 
  \sum\limits_{p<\infty} \,\sum\limits_{\ell\ge1} \, \hat\psi\left( {\ell\, \log p\over \log C_{\rho}
} \right)  {\log p\over \log C_{\rho}}
 	 \sum\limits_{k=1}^3 \frac{\ga_k(p)^\ell \, + \, \overline{\ga_k(p)^\ell}}{ p^{\ell/2}  }, &  \rho(\phi) = \phi,\\
& \\
B_{\rho(\phi)}
- 
 \sum\limits_{p<\infty} \,\sum\limits_{\ell\ge1}\, \hat\psi\left( {\ell \,\log p\over \log C_{\rho}
} \right)  {\log p\over \log C_{\rho}}\, 
 	\sum\limits_{1\le j \le k\le 3} \hskip -7pt \frac{\ga_j(p)^\ell{\ga_k(p)}^\ell \, + \, \overline{\ga_j(p)^\ell{\ga_k(p)}^\ell}}{ p^{\ell/2}  }, & \rho(\phi) = \sym^{2}\phi,\\ 
& \\
B_{\rho(\phi)}
- 
 2 \sum\limits_{p<\infty} \,\sum\limits_{\ell\ge1}\, \hat\psi\left( {\ell \,\log p\over \log C_{\rho}
} \right)  {\log p\over \log C_{\rho}} 
 	\sum\limits_{j=1}^3 \sum\limits_{k=1}^3 \frac{\ga_j(p)^\ell\,\overline{\ga_k(p)^\ell}}{ p^{\ell/2}  }, &  \rho(\phi) = \phi\times\bar\phi.
\end{cases}
\end{align}
with
\begin{equation}\label{B}
B_{\rho(\phi)} = 2r_{\rho(\phi)} \,\psi\left( {\log C_{\rho}\over 4\pi i}\right) +   {A_{\rho(\phi)} \over \log C_{\rho}
 } ,
 \end{equation}
 and where
\begin{align*}A_{\rho(\phi)} & =\hskip -2pt \begin{cases} -3\widehat\psi(0) \log\pi  +  \frac12\sum\limits_{k=1}^3  \int\limits_{-\infty}^\infty \bigg(  {\G'\over \G}\left( {\pi i x\over \log C_{\rho}%R
} + {1+2\alpha_k\over 4} \right)\\
&
\\
\hskip 160pt + {\G'\over \G}\left( {\pi i x\over \log C_{\rho}%R
} + {1-2\alpha_k\over 4} \right)\hskip -3pt\bigg)
\psi(x) \, dx, & \rho(\phi) = \phi,\\
& \\
 -5\widehat\psi(0) \log\pi  + \frac12 \sum\limits_{1\le j\le k\le 3}\;  \int\limits_{-\infty}^\infty \hskip -1pt\bigg(  {\G'\over \G}\left( {\pi i x\over \log C_{\rho}%R
} + {1+2(\alpha_j + \alpha_k)\over 4} \right)\\
&
\\ \hskip 145pt + {\G'\over \G}\left( {\pi i x\over \log C_{\rho}%R
} + {1-2(\alpha_j + \alpha_k)\over 4} \right)\hskip -3pt \bigg)
\psi(x) \, dx, & \rho(\phi) = \sym^{2}\phi,\\
 & \\ -9\widehat\psi(0) \log\pi  + \sum\limits_{j=1}^3 \sum\limits_{k=1}^3  \int\limits_{-\infty}^\infty {\G'\over \G}\left( {\pi i x\over \log C_{\rho}%R
} + {1+2(\alpha_j - \alpha_k)\over 4} \right)
\psi(x) \, dx, & \rho(\phi) = \phi\times \bar\phi. \end{cases}
\end{align*}

Now, for $\Re(\alpha) = 0,$ we have
$$\int\limits_{-\infty}^\infty \frac{\Gamma'}{\Gamma}\left(\frac{\pi ix}{\log C_{\rho}}  + \alpha + \frac14 \right) \psi(x) \, dx = \widehat\psi(0) \frac{\Gamma'}{\Gamma}\left(\alpha+\frac14\right) + \mathcal O\left( \Big(\big|\alpha+1/4\big|\, \log C_{\rho}\Big)^{-2}  \right). $$
If we combine this with the fact that 
$\frac{\Gamma'}{\Gamma}\left(\alpha + \frac14\right) = \log \alpha + \mathcal O(1)$
for $|\alpha| \ge \frac14$, 
it follows that
\begin{align*}A_{\rho(\phi)} & = \begin{cases}  \widehat\psi(0) \log\left( \pi^{-3}\cdot \underset {|\alpha_k| \,\ge \,\frac12} {\prod\limits_{1\le k\le 3}}\frac{|\alpha_k|}{2}\right) \; + \; \mathcal O(1), & \rho(\phi) = \phi,\\
& \\
 \widehat\psi(0) \log\left( \pi^{-5}\cdot \underset {|\alpha_j+\alpha_k| \,\ge \,\frac12} {\prod\limits_{1\le j\le k\le 3}}\frac{|\alpha_j+\alpha_k|}{2}\right) \; + \; \mathcal O(1), & \rho(\phi) = \sym^{2}\phi,\\
 & \\  \widehat\psi(0) \log\left( \pi^{-9}\cdot \underset {|\alpha_j-\alpha_k| \,\ge \,\frac12}{\prod\limits_{j=1}^3\prod\limits_{k=1}^3} \frac{|\alpha_j-\alpha_k|}{2}\right) \; + \; \mathcal O(1), & \rho(\phi) = \phi\times \bar\phi. \end{cases}
\end{align*}

Thus, in all cases, we have
\be\label{eq:ArhoPhi}
A_{\rho(\phi)} =  \widehat\psi(0) \log c_{\rho(\phi)} + \mathcal O(1)
,
\ee
where the analytic conductor $c_{\rho(\phi)}$ is given by \eqref{eq:AnalyCond}.

We combine the above formula for $A_{\rho(\phi)}$ with \eqref{Dphipsi} and \eqref{B}. The contribution to  \eqref{Dphipsi} from $\ell \ge 3$ is negligible (using the Ramanujan bound $|\ga_{j}(p)|\le1$), so we have
\begin{equation} \label{DphipsiRevised}
\boxed{D(\rho(\phi);\, \psi) = \widehat\psi(0) -\Sigma_{\rho(\phi)}^1 - \Sigma_{\rho(\phi)}^2 +  2r_{\rho(\phi)}\,\psi\left( {\log c_{\rho(\phi)}\over 4\pi i}\right) + \mathcal O\left({1\over \log c_{\rho(\phi)}}\right),}
\end{equation}
where $\Sigma_{\rho(\phi)}^\ell$ is the contribution from $\ell = 1,2,$ namely
  
 \begin{equation} \label{Sigma1}
\Sigma_{\rho(\phi)}^1 =\begin{cases}  - 
  \sum\limits_{p<\infty} \hat\psi\left( {\log p\over \log C_{\rho}
} \right)  {\log p\over \log C_{\rho}} 
 	 \sum\limits_{k=1}^3 \frac{\ga_k(p) \, + \, \overline{\ga_k(p)}}{ p^{1/2}  }, & {\rm if} \; \rho(\phi) = \phi,\\
& \\
- 
  \sum\limits_{p<\infty} \hat\psi\left( {\log p\over \log C_{\rho}
} \right)  {\log p\over \log C_{\rho}}\, 
 	\sum\limits_{1\le j \le k\le 3}  \frac{\ga_j(p) \ga_k(p) \, + \, \overline{\ga_j(p) \ga_k(p)}}{ p^{1/2}  }, & {\rm if}\; \rho(\phi) = \sym^{2}\phi,\\ 
& \\
- 
 2 \sum\limits_{p<\infty} \hat\psi\left( {\log p\over \log C_{\rho}
} \right)  {\log p\over \log C_{\rho}} 
 	\sum\limits_{j=1}^3 \sum\limits_{k=1}^3 \frac{\ga_j(p)\,\overline{\ga_k(p)}}{ p^{1/2}  }, & {\rm if} \; \rho(\phi) = \phi\times\bar\phi, \end{cases}
	\end{equation}
\vskip 4pt
\noindent
and

\begin{equation}\label{Sigma2}\Sigma_{\rho(\phi)}^2 =\begin{cases}  - 
  \sum\limits_{p<\infty} \hat\psi\left( {2\log p\over \log C_{\rho}
} \right)  {\log p\over \log C_{\rho}} 
 	 \sum\limits_{k=1}^3 \frac{\ga_k(p)^2 \, + \, \overline{\ga_k(p)^2}}{ p  }, & {\rm if} \; \rho(\phi) = \phi,\\
& \\
- 
  \sum\limits_{p<\infty} \hat\psi\left( {2\log p\over \log C_{\rho}
} \right)  {\log p\over \log C_{\rho}}\, 
 	\sum\limits_{1\le j \le k\le 3}  \frac{\ga_j(p)^2 \ga_k(p)^2 \, + \, \overline{\ga_j(p)^2 \ga_k(p)^2}}{ p  }, & {\rm if}\; \rho(\phi) = \sym^{2}\phi,\\ 
& \\
- 
 2 \sum\limits_{p<\infty} \hat\psi\left( {2\log p\over \log C_{\rho}
} \right)  {\log p\over \log C_{\rho}} 
 	\sum\limits_{j=1}^3 \sum\limits_{k=1}^3 \frac{\ga_j(p)^2 \,\overline{\ga_k(p)^2}}{ p  }, & {\rm if} \; \rho(\phi) = \phi\times\bar\phi.
\end{cases}
\end{equation}

%%%%%%%%%%%%%%%%%%%%%%%%%%%%%%%%%%%%%%%%
%%%%%%%%%%%%%%%%%%%%%%%%%%%%%%%%%%%%%%%%

\newpage

\section{Local Analysis}\label{subs:Loc}\
Recall the Fourier expansion of $\phi$:
 \be\label{eq:phiFourier}
\phi(z) = \sum_{\g\in U_2(\Z)\bk \SL_{2}(\Z)} \sum_{k_1\ge1}\sum_{k_{2}\neq0} 
{A_\phi(k_1,k_2)
\over 
 k_1 |k_2|
}\, 
W_{\nu}^{\sgn(k_{2})}\left(
\bp
k_{1}|k_{2}|&&\\
&k_{1}&\\
&&1
\ep
\bp
\g & \\
 & 1\\
\ep
z
\right)
,
 \ee
 where $\phi$ is normalized so that $A_\phi(1,1) = 1.$ Then the Fourier coefficients $A_\phi(k_1,k_2)$ satisfy the Hecke relations
 \be\label{hecke}
A_\phi(n,1) A_\phi(k_1,k_2) = \sum_{d_0d_1d_2=n\atop d_1|k_1,d_2|k_2} A_\phi\left({k_1d_0\over d_1},{k_2d_1\over d_2}\right),
\ee
as well as the conjugation relation $$\overline{A(k_1,k_2)} = A(k_2,k_1).$$
 Furthermore, the L-function associated to $\phi$ is given by
 $$L(s, \phi) = \prod_p L_p(s, \phi) = \,\sum_{n=1}^\infty \frac{A_\phi(1,n)}{n^s},$$
with  local factor
$$
L_p(s,\phi):= \prod_{j=1}^3\left(1-{\ga_j(p)\over p^s}\right)^{-1}
= \,\sum_{m\ge0}{A_\phi(1, p^m)\over p^{ms}} .
$$
It follows that 
$$
%\ga_1(p)+\ga_2(p)+\ga_3(p)=\gl(p)
A_\phi(1, p^m)=\sum_{u+v+w=m}\ga_1(p)^u\ga_2(p)^v\ga_3(p)^w.
$$
In particular,
$$
A_\phi(1,p)=\ga_1(p)+\ga_2(p)+\ga_3(p),
$$
and
\begin{align*}
A_\phi(1,p^2)
& =\ga_1(p)^2+\ga_2(p)^2+\ga_3(p)^2+\ga_1(p)\ga_2(p)+\ga_2(p)\ga_3(p)+\ga_1(p)\ga_3(p)\\
& =\ga_1(p)^2+\ga_2(p)^2+\ga_3(p)^2+\overline{A_\phi(1,p)}%\ga_1\ga_2+\ga_2\ga_3+\ga_1\ga_3
.
\end{align*}
Next, we analyze the local contribution at a prime $p$ occurring in the sum $\Sigma_{\rho(\phi)}^1$ given in \eqref{Sigma1}. There are 3 cases to consider. In the second and third case we apply the Hecke relations \eqref{hecke} to remove the product of two Fourier coefficients:
\begin{align*} \sum_{k=1}^3 \alpha_k(p) & = \boxed{A_\phi(1,p), }
\\ 
\sum_{1\le j\le k\le3} \alpha_j(p)\alpha_k(p) & = A_\phi(1,p)^2 - \alpha_1(p)\alpha_2(p) - \alpha_1(p)\alpha_3(p) - \alpha_2(p)\alpha_3(p)\\
& = A_\phi(1,p)^2 -\overline{A_\phi(1,p)} = A_\phi(1,p)^2 - A_\phi(p,1)\\
& = \boxed{A_\phi(1, p^2), }
\\
\sum_{j=1}^3\sum_{k=1}^3 \alpha_j(p) \overline{\alpha_k(p)} & = A_\phi(1,p) A_\phi(p,1) = \boxed{1 + A_\phi(p,p). }\end{align*}
Putting the above identities  into \eqref{Sigma1} yields

 \begin{equation} \label{Sigma1revised}
\Sigma_{\rho(\phi)}^1 =\begin{cases}  - 
  \sum\limits_{p<\infty} \hat\psi\left( {\log p\over \log C_{\rho}
} \right)  {\log p\over \log C_{\rho}} 
 	 \frac{A_\phi(1,p)\, + \, A_\phi(p,1)}{ p^{1/2}  }, & {\rm if} \; \rho(\phi) = \phi,\\
& \\
- 
  \sum\limits_{p<\infty} \hat\psi\left( {\log p\over \log C_{\rho}
} \right)  {\log p\over \log C_{\rho}}\, 
 	  \frac{A_\phi(1,p^2) \, + \, A_\phi(p^2,1)}{ p^{1/2}  }, & {\rm if}\; \rho(\phi) = \sym^{2}\phi,\\ 
& \\
- 
 2 \sum\limits_{p<\infty} \hat\psi\left( {\log p\over \log C_{\rho}
} \right)  {\log p\over \log C_{\rho}} 
 	 \frac{1 + A_\phi(p,p)}{ p^{1/2}  }, & {\rm if} \; \rho(\phi) = \phi\times\bar\phi, \end{cases}
	\end{equation}
\vskip 5pt
Next, we do the same for $\Sigma_{\rho(\phi)}^2$ given in \eqref{Sigma2}. As before, there are three cases to consider. We require the following Hecke relations obtained from \eqref{hecke}:
\beann
A_\phi(p,1)^2 
&=&
 A_\phi(p^2,1) + A_\phi(1,p),\\
A_\phi(1,p^2)^2
&=&
A_\phi(p^2,1) + A_\phi(p,p^2) + A_\phi(1,p^4)
,
\\
A_\phi(p,1) A_\phi(1,p^2)
&=&
A_\phi(p,p^2) + A_\phi(1, p^3),
\\
A_\phi(p^2,1)A_\phi(1,p^2)
&=&
A_\phi(1,1)+A_\phi(p,p)+A_\phi(p^2,p^2)
,
\\
-A_\phi(p^2,1)A_\phi(p,1)
&=&
-A_\phi(p,p)-A_\phi(p^3,1)
,
\\
-A_\phi(1,p^2)A_\phi(1,p)
&=&
-A_\phi(p,p)-A_\phi(1,p^3)
,
\\
A_\phi(p,1)A_\phi(1,p)
&=&
1+A_\phi(p,p)
.
\eeann
It follows from the above Hecke relations that:
\begin{align*}
\sum_{k=1}^3 \alpha_k(p)^2 & = \boxed{ A_\phi(1,p^2) - A_\phi(p,1), }\\
\sum_{1\le j\le k \le 3} \alpha_j(p)^2 \alpha_k(p)^2 & = \left(A_\phi(1,p^2) - A_\phi(p,1)\right)^2 - A_\phi(p^2,1) - A_\phi(p,1)\\
& = A_\phi(1,p^2)^2 - 2 A_\phi(p,1)A_\phi(1,p^2) + A_\phi(p,1)^2 - A_\phi(p^2,1) + A_\phi(p,1)
\\
& = \boxed{A_\phi(p^2,1) + A_\phi(1,p^4) -A_\phi(p,p^2) -2A_\phi(1,p^3) +A_\phi(1,p) + A_\phi(p,1), }
\\
\sum_{j=1}^3 \sum_{k=1}^3 \alpha_j(p)^2 \overline{\alpha_k(p)^2} & = \left(A_\phi(1,p^2) - A_\phi(p,1)\right)  \left(A_\phi(p^2, 1) - A_\phi(1,p)\right)\\
& = A_\phi(p^2,1) A_\phi(1,p^2) - A_\phi(p,1) A_\phi(p^2,1)\\
& \hskip 138pt - A_\phi(1,p^2) A_\phi(1,p) + A_\phi(p,1) A_\phi(1,p)
\\
& = \boxed{2 + A_\phi(p^2, p^2) - A_\phi(1,p^3) - A_\phi(p^3,1). }
\end{align*}
\vskip 5pt
Putting the above identities  into \eqref{Sigma2} yields

\begin{equation}\label{Sigma2revised}\Sigma_{\rho(\phi)}^2 =\begin{cases}  - 
  \sum\limits_{p<\infty} \hat\psi\left( {2\log p\over \log C_{\rho}
} \right)  {\log p\over \log C_{\rho}} 
 	\left( \frac{A_\phi(1,p^2) - A_\phi(p,1) \, + \, A_\phi(p^2,1) - A_\phi(1,p)}{ p  }\right), &  \rho(\phi) = \phi,\\
& \\
- 
  \sum\limits_{p<\infty}\hskip -2pt \hat\psi\left( {2\log p\over \log C_{\rho}
} \right) \hskip -2pt {\log p\over \log C_{\rho}}\hskip -2pt 
 	 \bigg( \frac{A_\phi(p^2,1) + A_\phi(1,p^4) -A_\phi(p,p^2) -2A_\phi(1,p^3) +A_\phi(1,p) - A_\phi(p,1)}{ p  }\\
	 &
	 \\ 
\hskip 90pt	 + \frac{A_\phi(1, p^2) + A_\phi(p^4, 1) -A_\phi(p^2, p) -2A_\phi(p^3, 1) +A_\phi(p, 1) - A_\phi(1, p)}{ p  }\bigg), &  \rho(\phi) = \sym^{2}\phi,\\ 
& \\
- 
  2\sum\limits_{p<\infty} \hat\psi\left( {2\log p\over \log C_{\rho}
} \right)  {\log p\over \log C_{\rho}} 
 	\left(\frac{2 + A_\phi(p^2, p^2) - A_\phi(1,p^3) - A_\phi(p^3,1)}{ p  }\right), &  \rho(\phi) = \phi\times\bar\phi.
\end{cases}
\end{equation}

\begin{remark}
{\rm 
It is much simpler to apply the Kuznetsov  formula to the sums $\Sigma_{\rho(\phi)}^1, \, \Sigma_{\rho(\phi)}^2,$ given
in \eqref{Sigma1revised} and \eqref{Sigma2revised} since the product of Fourier coefficients has been removed by the use of the Hecke relations. 
}
\end{remark}

%%%%%%%%%%%%%%%%%%%%%%%%%%%%%%%%%%%%%%%%
%%%%%%%%%%%%%%%%%%%%%%%%%%%%%%%%%%%%%%%%

\newpage

\section{Proof of Theorem \ref{thm:main}}\label{sec:pfMain}

%%%%%%%%%%%%%%%%%%%%%%%%%%%%%%%%%%%%%%%%
%%%%%%%%%%%%%%%%%%%%%%%%%%%%%%%%%%%%%%%%

%\newpage

%\subsection{Massaging $\gS_{1}$ and $\gS_{2}$}\label{subs:EstgS12}\

 By Parseval's Theorem,
$$
\int_\R\psi(x)W_{\rho(\phi)}(x)dx = \int_\R \hat\psi(y)\hat W_{\rho(\phi)}(y)dy,
$$
where $\hat\psi$ is the Fourier transform of $\psi$, and $\hat W_{\rho(\phi)}$ can be explicitly computed from \eqref{Wis} as a distribution:
$$
\hat W_{\rho(\phi)}(y)=\begin{cases} \delta_0(y), & {\rm if } \;\rho(\phi)=\phi \; {\rm or} \; \sym^{2}\phi ,\\
  \gd_0(y)-\scriptstyle{\threecase{\scriptstyle1/2,}{$\scriptstyle |y|<1,$}{\scriptstyle 1/4,}{$\scriptstyle |y|=\pm1,$}{\scriptstyle 0,}{ $\scriptstyle |y|>1,$.}} & {\rm if} \; \rho(\phi) =\Ad\phi.\end{cases}
$$
As the support of $\hat\psi$ will be restricted well inside $(-1,1)$, it follows that
\begin{equation}
\int_\R \hat\psi(y)\hat W_{\rho(\phi)}(y)dy = 
\begin{cases} 
%\int_\R \hat\psi(y)\gd_0(y)dy =
 \widehat\psi(0), & {\rm if} \; \rho(\phi) = \phi,\\
%\int_\R \hat\psi(y)\gd_0(y)dy =
 \widehat\psi(0), & {\rm if} \; \rho(\phi) = \sym^2\phi\\
%\int_\R \hat\psi(y)\gd_0(y)dy - \int_{supp(\hat\psi)}\hat\psi(y) \foh dy,  \\
%\hskip 75pt=
  \hat\psi(0) - \foh \psi(0), &  {\rm if} \; \rho(\phi) = \Ad\phi.
\end{cases}
\label{Wtransform}\end{equation}

\vskip 10pt
Recall the asymptotic formula \eqref{DphipsiRevised} for the low-lying zeros sum:
\begin{equation} \label{DphipsiRevised2}
\boxed{D(\rho(\phi); \,\psi) = \widehat\psi(0) -\Sigma_{\rho(\phi)}^1 - \Sigma_{\rho(\phi)}^2 +  2r_{\rho(\phi)}\psi\left( {\log c_{\rho(\phi)}\over 4\pi i}\right) + \mathcal O\left({1\over \log c_{\rho(\phi)}}\right),}
\end{equation}
where by \eqref{residues}, we have

\begin{equation*} r_{\rho(\phi)} = \begin{cases} 0 & {\rm if} \;  \rho(\phi) = \phi,
\\
 0 & {\rm if} \;  \rho(\phi) = \sym^{2}\phi \; {\rm and} \; \phi \; {\rm not \; self \; dual,}\\
 1 & {\rm if} \; \rho(\phi) = \phi\times\bar\phi, \end{cases}
 \end{equation*}
 and $\Sigma_{\rho(\phi)}^1$, $\Sigma_{\rho(\phi)}^2$ are given by \eqref{Sigma1revised}, and 
\eqref{Sigma2revised}, respectively.

We will prove Theorem \ref{thm:main} where the limiting density $\int_\R \hat\psi(y)\hat W(y)dy$ is in  the form
\eqref{Wtransform}. Note that $\widehat\psi(0)$ already appears in \eqref{DphipsiRevised2}.

\vskip 10pt\noindent
$\underline{\text{\bf The case when $\rho(\phi) = \phi\times\bar\phi$\,:}}  $
\vskip 10pt
The main contribution to $\Sigma_{\phi\times\bar\phi}^1 + \Sigma_{\phi\times\bar\phi}^2$ 
%will
 comes from
\be\label{twosteps}
2 \sum_{p%<\infty
 } 
 \frac 1{\sqrt p}\;
\hat\psi\left( { \log p\over \log C_{\rho}
} \right)  {\log p\over \log C_{\rho}
}
\; + \; 4 \sum_{p%<\infty
 } \frac{1}{p} \;
 \hat\psi\left( {2 \log p\over \log C_{\rho}
} \right)  {\log p\over\log C_{\rho}
}.
\ee
This may be computed in two steps. 
\vskip 10pt
\noindent
{\bf Step 1:} We apply the explicit formula \eqref{eq1}  to the function $$\Lambda(s) = \pi^{-\frac{s}{2}} \Gamma\left(\frac{s}{2}\right) \zeta(s).$$  Note that for $\log C_{\rho}$ sufficiently large, we have$$\sum_{\Lambda(\rho) = 0} G(\rho) = \sum_{\Lambda(\frac12+i\gamma)=0} \psi\left(\frac{\gamma \,\log C_{\rho}}{2\pi}\right) = 0,$$
since there are no low-lying (with $|\gamma| \ll 1/\log C_{\rho}$) non-trivial zeros of the Riemann zeta function.
Further $G(0)+G(1) = 2\psi\left(\frac{\log C_{\rho}}{4\pi i}\right)$. It follows that the explicit formula for the Riemann zeta function takes the form:
\begin{align*}  2\psi\left(\frac{\log C_{\rho}}{4\pi i}\right) & = 2\sum_{\ell\ge1}\frac{1}{p^{\ell/2}} \cdot\widehat\psi\left(\frac{\ell\, \log p}{\log C_{\rho}} \right)\frac{\log p}{\log C_{\rho}} + \frac{\widehat\psi(0)}{\log C_{\rho}}\frac{\log\pi}{2}\\
&\hskip 90pt  - \frac{1}{\log C_{\rho}} \int_{-\infty}^\infty \frac{\Gamma'}{\Gamma} \left( \frac{2\pi i x}{\log C_{\rho}} + \frac{1}{4}  \right) \psi(x) \; dx\\
&
\\& = 2\sum_{\ell=1}^2\frac{1}{p^{\ell/2}} \cdot\widehat\psi\left(\frac{\ell\,\log p}{\log C_{\rho}} \right)\frac{\log p}{\log C_{\rho}} \; + \; \mathcal O\left(\frac{1}{\log C_{\rho}}  \right).\end{align*}
Note that the above accounts for a large part of \eqref{twosteps}
 and almost cancels the term $2\psi\left(\frac{\log c_{\rho(\phi)}}{4\pi i}\right)$ in \eqref{DphipsiRevised2}. In fact, we get
 \begin{align}
 D(\phi\times\bar\phi; \,\psi) & = \widehat\psi(0)  - \Sigma_{\phi\times\bar\phi}^2 +  2\psi\left( {\log c_{\phi\times\bar\phi}\over 4\pi i}\right) - 2\psi\left( {\log C_{\rho}\over 4\pi i}\right)\notag\\
 &
 \hskip 38pt - 2 \sum_{p
 } \frac{1}{p}\;
 \hat\psi\left( {2 \log p\over \log C_{\rho}
} \right)  {\log p\over\log C_{\rho}
} + \mathcal O\left({1\over \log c_{\rho(\phi)}}\right).\label{DphipsiRevised3}\end{align}

 What is really happening here is that every Rankin-Selberg $L$-function $L(s, \phi\times\bar\phi)$ in the family is divisible by the Riemann zeta function which has a pole at $s = 1.$ But the Riemann zeta function does not contribute low lying zeros, so the contribution from the pole at $s =1$ is cancelled.

\vskip 10pt\noindent
{\bf Step 2:} In the second step we make use of the classical Riemann hypothesis which implies that $\Psi(x) = x+\mathcal O\left(x^{\frac12+\epsilon}\right).$ It follows that\begin{align*} 
& 2 \sum_{p} \frac{1}{p}\;
 \hat\psi\left( {2 \log p\over \log C_{\rho}} \right)  {\log p\over \log C_{\rho}}
 =
\frac{2}{\log C_{\rho}} \int_1^\infty
 \hat\psi\left( \frac{2 \log x}{\log C_{\rho}} \right) x^{-1} d\Psi(x)
\\ &\\
& \hskip 30pt =
{2\over \log C_{\rho}%R
} \int_1^\infty
 \hat\psi\left( {2 \log x\over \log C_{\rho}%R
}\, \right) x^{-2} \left(x+\mathcal O\left(x^{\frac12+\epsilon}\right)\right)\, dx\\
& \hskip 53pt
-
{2\over \log C_{\rho}%R
} \int_1^\infty
{2 \over x\ \log C_{\rho}%R
}\, 
 \hat\psi'\left( {2 \log x\over \log C_{\rho}%R
} \right) x^{-1} \left(x+\mathcal O\left(x^{\frac12+\epsilon}\right)\right)\, dx
\\
& \\
&\hskip 30pt =
\int_0^\infty
 \hat\psi\left(u  \right) du
\; + \; \mathcal O\left(\frac{1}{\log C_{\rho}}   \right)
\\ 
&\hskip 30pt =
\frac12 \psi(0)
\; + \; \mathcal O\left(\frac{1}{\log C_{\rho}}   \right).
\end{align*}

\vskip 10pt
If we now combine the computation in Step 2 with equation \eqref{DphipsiRevised3},  it follows that

\be\label{phitimesbarphi}
\boxed{
D(\phi\times\bar\phi; \,\psi) = \hat\psi(0)-\frac12 \psi(0) + 2\psi\left(\frac{\log c_{\phi\times\bar\phi}}{4\pi i}\right) - 2\psi\left(\frac{\log C_{\rho}}{4\pi i}\right) +\Sigma_{\phi\times\bar\phi}^3
+ O\left({1\over \log c_{\phi\times\bar\phi}}\right),}
\ee

where
\begin{align} \label{Sigma3} \Sigma_{\phi\times\bar\phi}^3 & = 2 \sum_p \Bigg[\frac{A_\phi(p, p)}{\sqrt{p}} \,\hat\psi\left(\frac{\log p}{\log C_{\rho}}\right)\nonumber\\
& \hskip 30pt + \left(\frac{A_\phi(p^2, p^2)-A_\phi(1, p^3) - A_\phi(p^3, 1)}{p}\right) \hat\psi\left(\frac{2\log p}{\log C_{\rho}}\right)\Bigg]  \frac{\log p}{\log C_{\rho}}.
\end{align}
\vskip 10pt
To prove Theorem \ref{thm:main} for the family of Rankin Selberg L-functions, we make use of \eqref{phitimesbarphi} and the decay properties of $h_{T,R}$  to obtain
$$\sum_{j}D(\phi_j\times\phi_j; \,\psi)\, {h_{T,R}(\nu^{(j)})\over \cL_{j}} = \sum_j\Bigg[\hat\psi(0)-\frac12 \psi(0) + \Sigma_{\phi_j\times\bar\phi_j}^3
+ O\left({\log\log T\over \log T}\right)\Bigg]  {h_{T,R}(\nu^{(j)})\over \cL_{j}}.$$
as we average over Maass forms $\phi_j\; (j=1,2,\ldots)$. The term $\log\log T/\log T$ arises after breaking the sum into two pieces corresponding to  $\phi_j$ with conductor $c_{\phi_j} \ll T^3/\log T$ and $c_{\phi_j} \gg T^3/\log T$. To evaluate the above sum, it remains to estimate 
$$ \sum_j \Sigma_{\phi_j\times\bar\phi_j}^3 {h_{T,R}(\nu^{(j)})\over \cL_{j}}.$$ 
Say the support of $\hat\psi$ is in $(-\gd,\gd)$. It immediately follows from \eqref{Sigma3} that 
\begin{align*} &  \sum_j \Sigma_{\phi_j\times\bar\phi_j}^3 {h_{T,R}(\nu^{(j)})\over \cL_{j}} =  2\sum_j \Bigg[\sum_{p\ll T^{6 \delta}} \frac{A_\phi(p, p)}{\sqrt{p}} \hat\psi\left(\frac{\log p}{\log C_{\rho}}\right)  \frac{\log p}{\log C_{\rho}}\\
& \hskip 15pt + \sum_{p \ll T^{3 \delta}}\left(\frac{A_\phi(p^2, p^2)-A_\phi(1, p^3) - A_\phi(p^3, 1)}{p}\right) \hat\psi\left(\frac{2\log p}{\log C_{\rho}}\right)  \frac{\log p}{\log C_{\rho}} \Bigg] \cdot   {h_{T,R}(\nu^{(j)})\over \cL_{j}}.\end{align*}
To finish the estimation, we apply Theorem \ref{thm:Kuz} to obtain
\begin{align*}\sum_j \Sigma_{\phi_j\times\bar\phi_j}^3 {h_{T,R}(\nu^{(j)})\over \cL_{j}}
&  \ll 
 T^{3+3R+\vep} \left[\sum_{p\ll T^{6 \delta}} p^{4-\frac12}\;   \frac{\log p}{\log C_{\rho}} \;+\;
 \sum_{p \ll T^{3 \delta}} p^7 \,  \frac{\log p}{\log C_{\rho}} \right] \\
 & \ll T^{3 + 3R + 27\delta +\vep}.
 \end{align*}
 So we need
 $$\delta < \frac2{27}.$$
 With this choice of $\delta$ we obtain
 $$\sum_{j}D(\phi_j\times\phi_j; \,\psi)\, {h_{T,R}(\nu^{(j)})\over \cL_{j}} = \Bigg[\hat\psi(0)-\frac12 \psi(0) + O\left({\log\log T\over \log T}\right)\Bigg]\cdot \sum_j {h_{T,R}(\nu^{(j)})\over \cL_{j}},$$
 as claimed.
 
 \begin{remark}
{\rm 
The family  of Rankin-Selberg L-functions for $GL(3)$ has the same symmetry type as the family of adjoint L-functions in view of the identity
$L(s, \Ad\phi) = \frac{L(s,\; \phi\times \bar\phi)}{\zeta(s)},$
and the fact that $\zeta(s)$ has no low-lying zeros.}
\end{remark}

 \vskip 10pt\noindent
$\underline{\text{\bf The case when $\rho(\phi) =\sym^{2}\phi$\,:}}  $
\vskip 10pt

It immediately follows from \eqref{DphipsiRevised2}, \eqref{Sigma1revised}, and \eqref{Sigma2revised}, that
$$
D(\sym^{2}\phi; \,\psi) = \widehat\psi(0) -\Sigma_{\sym^{2}\phi}^1 - \Sigma_{\sym^{2}\phi}^2 +  2r_{\sym^{2}\phi} \,\psi\left( {\log c_{\sym^{2}\phi}\over 4\pi i}\right) + \mathcal O\left({1\over \log c_{\sym^{2}\phi}}\right),
$$
where
$$\Sigma_{\sym^{2}\phi}^1 = 
 - \sum\limits_{p<\infty} \hat\psi\left( {\log p\over \log C_{\rho}
} \right)  {\log p\over \log C_{\rho}}\, 
 	  \frac{A_\phi(1,p^2) + A_\phi(p^2,1)}{ p^{1/2}  }$$
and
\begin{align*}\Sigma_{\sym^{2}\phi}^2  & \; = \; - \hskip -2pt
  \sum\limits_{p<\infty}\hskip -2pt \hat\psi\left( {2\log p\over \log C_{\rho}
} \right)  {\log p\over \log C_{\rho}}\\
& \hskip 35pt\times 
 	  \bigg(\frac{A_\phi(p^2,1) + A_\phi(1,p^4) -A_\phi(p,p^2) -2A_\phi(1,p^3) +A_\phi(1,p) - A_\phi(p,1)}{ p  }\\
	  &
	\hskip 66pt  + \frac{A_\phi(1, p^2) + A_\phi(p^4,1) -A_\phi(p^2,p) -2A_\phi(p^3,1) +A_\phi(p,1) - A_\phi(1,p)}{ p  }\bigg).
	  \end{align*}
	  In this case $r_{\sym^{2}\phi} = 0$ unless $\phi$ is self dual, which happens only if $\phi$ is a symmetric square lift from $GL(2)$. This occurs for $\asymp T^2$ cases out of $\asymp T^6$, and hence contributes a negligible error term to the low-lying zeros sum. The method to estimate  $\Sigma_{\sym^{2}\phi}^1$ and $\Sigma_{\sym^{2}\phi}^2$, using Theorem \ref{thm:Kuz}, is exactly the same as the method used above for the case of $\rho(\phi) = \phi\times\bar\phi.$ Because of the presence of $A_\phi(1,p^2), A_\phi(p^2, 1)$ in $\Sigma_{\sym^{2}\phi}^1$ and $A_\phi(1,p^4), A_\phi(p^4, 1)$ in $\Sigma_{\sym^{2}\phi}^2$, we obtain the same value $\delta < \frac2{27}$ (which we obtained for the case  $\rho(\phi) = \phi\times\bar\phi$) for the support of $\widehat\psi.$ It follows that
$$\sum_{j}D(\sym^{2}\phi; \,\psi)\, {h_{T,R}(\nu^{(j)})\over \cL_{j}} = \Bigg[\hat\psi(0)
+ O\left({\log\log T\over \log T}\right)\Bigg] \cdot \sum_j {h_{T,R}(\nu^{(j)})\over \cL_{j}}.$$

 \vskip 10pt\noindent
$\underline{\text{\bf The case when $\rho(\phi) =\phi$\,:}}  $
\vskip 10pt
In this case, we have that the residue $r_{\phi} = 0.$ It then follows from \eqref{DphipsiRevised2}, \eqref{Sigma1revised}, and \eqref{Sigma2revised}, that
$$
D(\phi; \,\psi)  = \widehat\psi(0) -\Sigma_{\phi}^1 - \Sigma_{\phi}^2 +    \mathcal O\left({1\over \log c_{\phi}}\right)$$
where $$\Sigma_{\phi}^1 =   -2 \sum\limits_{p<\infty} \hat\psi\left( {\log p\over \log C_{\rho}
} \right)  {\log p\over \log C_{\rho}} 
 	 \frac{A_\phi(1,p) + A_\phi(p,1)}{ p^{1/2}  }$$
	 and
$$\Sigma_{\phi}^2 = - 
  \sum\limits_{p<\infty} \hat\psi\left( {2\log p\over \log C_{\rho}
} \right)  {2\log p\over \log C_{\rho}} \;
 	 \frac{A_\phi(1,p^2) - A_\phi(p,1) \, + \, A_\phi(p^2,1) - A_\phi(1,p)}{ p  }.$$
	 Assume the support of $\hat\psi$ is in $(-\gd,\gd)$. As before, we have
	 $$\sum_{j}D(\phi_j; \,\psi)\, {h_{T,R}(\nu^{(j)})\over \cL_{j}} = \sum_j\Bigg[\hat\psi(0) - \Sigma_{\phi_j}^1 -  \Sigma_{\phi_j}^2
+ O\left({\log\log T\over \log T}\right)\Bigg]  {h_{T,R}(\nu^{(j)})\over \cL_{j}},$$	 
and invoking  Theorem \ref{thm:Kuz} again gives
\begin{align*} & \sum_j  \left(-\Sigma_{\phi_j}^1 - \Sigma_{\phi_j}^2\right)  \,{h_{T,R}(\nu^{(j)})\over \cL_{j}}  =  \sum_j \Bigg[ \sum\limits_{p\ll T^{3\delta}} \hat\psi\left( {\log p\over \log C_{\rho}
} \right)  {\log p\over \log C_{\rho}} 
 	 \frac{A_\phi(1,p) + A_\phi(p,1)}{ p^{1/2}  }\\
	 &
	 \hskip 20pt + \sum\limits_{p\ll T^{3\delta/2 }} \hat\psi\left( {2\log p\over \log C_{\rho}
} \right)  {\log p\over \log C_{\rho}} \;
 	 \frac{A_\phi(1,p^2) - A_\phi(p,1) \, + \, A_\phi(p^2, 1) - A_\phi(1,p)}{ p  }\Bigg] \,{h_{T,R}(\nu^{(j)})\over \cL_{j}}\\
	 & \\
	 &
\hskip 100pt \ll	 T^{3+3R+\vep}\left[\sum_{p \ll T^{3\delta}} p^{3/2}  + \sum_{p\ll T^{3\delta/2}} p^3  \right]\\
&\\
& \hskip 100pt\ll T^{3+3R+\vep}\, T^{\frac{15}{2}\delta}.
\end{align*} 
So in this case, we need $\delta < \frac{4}{15}.$

This completes the proof of Theorem \ref{thm:main}.

%%%%%%%%%%%%%%%%%%%%%%%%%%%%%%%%%%%%%%%%
%%%%%%%%%%%%%%%%%%%%%%%%%%%%%%%%%%%%%%%%

%%%%%%%%%%%%%%%%%%%%%%%%%%%%%%%%%%%%%%%%
%%%%%%%%%%%%%%%%%%%%%%%%%%%%%%%%%%%%%%%%

%\newpage

\bibliographystyle{alpha}

\bibliography{AKbibliog}

\newcommand{\etalchar}[1]{$^{#1}$}
\begin{thebibliography}{AIL{\etalchar{+}}11}

\bibitem[AIL{\etalchar{+}}11]{MillerEtAl2011}
N.~Amersi, G.~Iyer, O.~Lazarev, S.~Miller, and L.~Zhang.
\newblock Low-lying zeros of cuspidal {M}aass forms, 2011.
\newblock Preprint. arXiv:1111.6524v2.

\bibitem[BF90]{BumpFriedberg1990}
Daniel Bump and Solomon Friedberg.
\newblock The exterior square automorphic {$L$}-functions on {${\rm GL}(n)$}.
\newblock In {\em Festschrift in honor of I. I. Piatetski-Shapiro on the
  occasion of his sixtieth birthday, Part II (Ramat Aviv, 1989)}, volume~3 of
  {\em Israel Math. Conf. Proc.}, pages 47--65. Weizmann, Jerusalem, 1990.

\bibitem[BFG88]{BumpFriedbergGoldfeld1988}
Daniel Bump, Solomon Friedberg, and Dorian Goldfeld.
\newblock Poincar\'e series and {K}loosterman sums for {${\rm SL}(3,{\bf Z})$}.
\newblock {\em Acta Arith.}, 50(1):31--89, 1988.

\bibitem[BG92]{BumpGinzburg1992}
D.~Bump and D.~Ginzburg.
\newblock Symmetric square ${L}$-functions on {G}{L}$(r)$.
\newblock {\em Ann. of Math. (2)}, 136(1):137--205, 1992.

\bibitem[Blo11]{Blomer2011}
V.~Blomer.
\newblock Applications of the {K}uznetsov formula on ${GL}(3)$, 2011.
\newblock Preprint.

\bibitem[DM06]{DuenezMiller2006}
Eduardo Due{\~n}ez and Steven~J. Miller.
\newblock The low lying zeros of a {$\rm GL(4)$} and a {$\rm GL(6)$} family of
  {$L$}-functions.
\newblock {\em Compos. Math.}, 142(6):1403--1425, 2006.

\bibitem[DM09]{DuenezMiller2009}
Eduardo Due{\~n}ez and Steven~J. Miller.
\newblock The effect of convolving families of {$L$}-functions on the
  underlying group symmetries.
\newblock {\em Proc. Lond. Math. Soc. (3)}, 99(3):787--820, 2009.

\bibitem[GJ72]{GodementJacquet1972}
Roger Godement and Herv{\'e} Jacquet.
\newblock {\em Zeta functions of simple algebras}.
\newblock Lecture Notes in Mathematics, Vol. 260. Springer-Verlag, Berlin,
  1972.

\bibitem[GK11]{GoldfeldKontorovich2011}
D.~Goldfeld and A.~Kontorovich.
\newblock On the determination of the {P}lancherel measure for
  {L}ebedev-{W}hittaker transforms on {GL}($n$), 2011.
\newblock To appear, {\it Acta Arithmetica}.

\bibitem[GLS04]{GelbartLapidSarnak2004}
S.~Gelbart, E.~Lapid, and P.~Sarnak.
\newblock A new method for lower bounds of ${L}$-functions.
\newblock {\em C. R. Math. Acad. Sci. Paris}, 339(2):91--94, 2004.

\bibitem[Gol06]{Goldfeld2006}
Dorian Goldfeld.
\newblock {\em Automorphic forms and {$L$}-functions for the group {${\rm
  GL}(n,\bold R)$}}, volume~99 of {\em Cambridge Studies in Advanced
  Mathematics}.
\newblock Cambridge University Press, Cambridge, 2006.
\newblock With an appendix by Kevin A. Broughan.

\bibitem[Gul05]{Gulolu2005}
A.~Gulolu.
\newblock Low-lying zeroes of symmetric power ${L}$-functions.
\newblock {\em IMRN}, (9):517--550, 2005.

\bibitem[HL94]{HoffsteinLockhart1994}
Jeffrey Hoffstein and Paul Lockhart.
\newblock Coefficients of {M}aass forms and the {S}iegel zero.
\newblock {\em Ann. of Math. (2)}, 140(1):161--181, 1994.
\newblock With an appendix by Dorian Goldfeld, Hoffstein and Daniel Lieman.

\bibitem[HM07]{HughesMiller2007}
C.~P. Hughes and Steven~J. Miller.
\newblock Low-lying zeros of {$L$}-functions with orthogonal symmetry.
\newblock {\em Duke Math. J.}, 136(1):115--172, 2007.

\bibitem[HR95]{HoffsteinRamakrishnan1995}
Jeffrey Hoffstein and Dinakar Ramakrishnan.
\newblock Siegel zeros and cusp forms.
\newblock {\em Internat. Math. Res. Notices}, 6:279--308, 1995.

\bibitem[IK04]{IwaniecKowalski}
Henryk Iwaniec and Emmanuel Kowalski.
\newblock {\em Analytic number theory}, volume~53 of {\em American Mathematical
  Society Colloquium Publications}.
\newblock American Mathematical Society, Providence, RI, 2004.

\bibitem[ILS00]{IwaniecLuoSarnak2000}
Henryk Iwaniec, Wenzhi Luo, and Peter Sarnak.
\newblock Low lying zeros of families of {$L$}-functions.
\newblock {\em Inst. Hautes \'Etudes Sci. Publ. Math.}, (91):55--131 (2001),
  2000.

\bibitem[JPSS83]{JacquetPSShalika1983}
H.~Jacquet, I.~Piatetskii-Shapiro, and J.~Shalika.
\newblock Rankin-{S}elberg convolutions.
\newblock {\em Amer. J. Math}, 105(2):367--464, 1983.

\bibitem[JS90]{JacquetShalika1990}
Herv{\'e} Jacquet and Joseph Shalika.
\newblock Exterior square {$L$}-functions.
\newblock In {\em Automorphic forms, Shimura varieties, and $L$-functions,
  Vol.\ II (Ann Arbor, MI, 1988)}, volume~11 of {\em Perspect. Math.}, pages
  143--226. Academic Press, Boston, MA, 1990.

\bibitem[JS77]{JacquetShalika1976}
H.~Jacquet and J.~Shalika.
\newblock A non-vanishing theorem for zeta functions of {GL}$(n)$.
\newblock {\em Invent. Math}, 38(1):1--16, 1976/77.

\bibitem[Kon10]{Kontorovich2010}
Alex~V. Kontorovich.
\newblock The {D}irichlet series for the exterior square {$L$}-function on
  {${\rm GL}(n)$}.
\newblock {\em Ramanujan J.}, 21(3):263--266, 2010.

\bibitem[KS99]{KatzSarnak1999}
Nicholas~M. Katz and Peter Sarnak.
\newblock {\em Random matrices, {F}robenius eigenvalues, and monodromy},
  volume~45 of {\em American Mathematical Society Colloquium Publications}.
\newblock American Mathematical Society, Providence, RI, 1999.

\bibitem[Luo01]{Luo2001}
Wenzhi Luo.
\newblock Nonvanishing of {$L$}-values and the {W}eyl law.
\newblock {\em Ann. of Math. (2)}, 154(2):477--502, 2001.

\bibitem[Mil01]{Miller2001}
Stephen~D. Miller.
\newblock On the existence and temperedness of cusp forms for {${\rm SL}\sb
  3({\Bbb Z})$}.
\newblock {\em J. Reine Angew. Math.}, 533:127--169, 2001.

\bibitem[P{\v{S}}75]{PS1975}
I.~I. Pjateckij-{\v{S}}apiro.
\newblock Euler subgroups.
\newblock In {\em Lie groups and their representations ({P}roc. {S}ummer
  {S}chool, {B}olyai {J}\'anos {M}ath. {S}oc., {B}udapest, 1971)}, pages
  597--620. Halsted, New York, 1975.

\bibitem[Roy01]{Royer2001}
E.~Royer.
\newblock Petits zeros de fonctions ${L}$ de formes modulaires.
\newblock {\em Acta Arith}, 99(2):147--172, 2001.

\bibitem[Sar08]{Sarnak2008}
P.~Sarnak.
\newblock Definition of families of ${L}$-functions, 2008.
\newblock Online note, available at {\it
  http://publications.ias.edu/sarnak/paper/507}.

\bibitem[Sha73]{Shalika1973}
J.~A. Shalika.
\newblock On the multiplicity of the spectrum of the space of cusp forms of
  {${\rm GL}_{n}$}.
\newblock {\em Bull. Amer. Math. Soc.}, 79:454--461, 1973.

\bibitem[ST12]{ShinTemplier2012}
S.~Shin and N.~Templier, 2012.
\newblock Preprint.

\bibitem[Sta01]{Stade2001}
Eric Stade.
\newblock Mellin transforms of {${\rm GL}(n,\Bbb R)$} {W}hittaker functions.
\newblock {\em Amer. J. Math.}, 123(1):121--161, 2001.

\bibitem[Sta02]{Stade2002}
Eric Stade.
\newblock Archimedean {$L$}-factors on {${\rm GL}(n)\times{\rm GL}(n)$} and
  generalized {B}arnes integrals.
\newblock {\em Israel J. Math.}, 127:201--219, 2002.

\bibitem[Wal92]{Wallach1992}
Nolan~R. Wallach.
\newblock {\em Real reductive groups. {II}}, volume 132 of {\em Pure and
  Applied Mathematics}.
\newblock Academic Press Inc., Boston, MA, 1992.

\bibitem[You06]{Young2006}
M.~Young.
\newblock Low-lying zeros of families of elliptic curves.
\newblock {\em JAMS}, 19(1):205--250, 2006.

\end{thebibliography}

\end{document}